\documentclass[12pt,reqno]{amsart}
\usepackage{amsthm,amsfonts,amssymb,euscript}

\theoremstyle{plain}
  \newtheorem{theorem}[subsection]{Theorem}
  
  \newtheorem{proposition}[subsection]{Proposition}
  \newtheorem{lemma}[subsection]{Lemma}
  \newtheorem{corollary}[subsection]{Corollary}
\def\RR{\mathbb R}
\def\bphi{{\boldsymbol{\phi} }}
\def\bpsi{{\boldsymbol{\psi} }}
\def\F{{\mathcal F}}
\def\R{{\mathbb{R}}}

\def\A{{\mathcal A}}
\def\B{{\mathcal B}}

\def\Q{{\mathcal Q}}
\def\H{{\mathcal H}}
\def\N{{\mathcal N}}

\def\E{{\mathcal E}}

\def\ch{\mbox{ch} (k)}
\def\chbb{\mbox{chb} (k)}
\def\trigh{\mbox{trigh} (k)}

\def\sh{\mbox{sh} (k)}
\def\shh{\mbox{sh}}
\def\chh{\mbox{ch}}
\def\shhb{\overline {\mbox{sh}}}
\def\th{\mbox{th} (k)}
\def\chhb{\overline {\mbox{  ch}}}
\def\shbb{\mbox{shb} (k)}

\def\shb{\overline{\mbox{sh}(k)}}
\def\z{\zeta}
\def\zb{\overline{\zeta}}
\def\FF{{\mathbb{F}}}
\def\L{\Lambda}

\theoremstyle{remark}
  \newtheorem{remark}[subsection]{Remark}
  
\theoremstyle{definition}

\newcommand{\llp}[1]{ \Vert \, #1 \, \Vert }
\newcommand{\rr}{\mathbb{R}}

\newcommand{\bd}{\partial}
\newcommand{\lapl}{\Delta}

\newcommand{\tr}{\operatorname{tr}}
\newcommand{\grad}{\nabla}
\newcommand{\vect}[1]{\mathbf{ #1 }}

\newcommand{\I}{\text{i}}

\begin{document}

\def\vac{{\big\vert 0\big>}}
\def\fpsi{{\big\vert\psi\big>}}
\def\bphi{{\boldsymbol{\phi} }}
\def\bpsi{{\boldsymbol{\psi} }}
\def\F{{\mathcal F}}
\def\R{{\mathbb{R}}}
\def\I{{\mathcal I}}
\def\Q{{\mathcal Q}}
\def\ch{\mbox{\rm ch} (k)}
\def\cht{\mbox{\rm ch} (2k)}
\def\chbb{\mbox{\rm chb} (k)}
\def\trigh{\mbox{\rm trigh} (k)}
\def\p{\mbox{p} (k)}
\def\sh{\mbox{\rm sh} (k)}
\def\sht{\mbox{\rm sh} (2k)}
\def\shh{\mbox{\rm sh}}
\def\chh{\mbox{ \rm ch}}
\def\th{\mbox{\rm th} (k)}
\def\thb{\overline{\mbox{\rm th} (k)}}
\def\shhb{\overline {\mbox{\rm sh}}}
\def\chhb{\overline {\mbox{\rm ch}}}
\def\shbb{\mbox{\rm shb} (k)}
\def\chbt{\overline{\mbox{\rm ch}  (2k)}}
\def\shb{\overline{\mbox{\rm sh}(k)}}
\def\shbt{\overline{\mbox{\rm sh}(2k)}}
\def\z{\zeta}
\def\zb{\overline{\zeta}}
\def\FF{{\mathbb{F}}}
\def\S{{\bf S}}
\def\Sr{{\bf S}_{red}}
\def\W{{\bf W}}
\def\WW{{\mathcal W}}
\def\L{{\mathcal L}}
\def\E{{\mathcal E}}
\def\X{{\tilde X}}
\def\Nor{{\rm Nor}}

\def \kb{{\bf k}}
\def\x{{\bf x}}
\def\fhat{{\widehat{f}}}
\def \ghat{{\widehat{g}}}
\def \fhatb{{\widehat{f}^{\ast}}}
\def \ghatb{{\widehat{g}^{\ast}}}
\def \Sb{\mathbb S}
\def\Zb {\mathbb Z}
\def\Lm{\Lambda}
\def \Lmb{\overline{\Lambda}}
\def \Gam{\Gamma}
\def\Gamb{\overline{\Gamma}}
\def\Lmh{\widehat{\Lambda}}
\def\Fhat{\widehat{F}}
\def\Fhatb{\overline{\widehat{F}}}
\def\E{{\mathcal E}}
\def\Sop{{\bf S}^{t}_{x_{1},x_{2}}}
\def\Wop{{\bf S}^{\pm,t}_{x_{1},x_{2}}}
\def\Sp{{\bf S}^{t}_{x_{1}}}
\def\nb{\nabla}
\def\Dl{\Delta}
\def\Acal{{\mathcal A}}
\def\ub{\overline{u}}
\def\tr{{\rm tr}}
\def\dig{{\rm diag}}

\def\Hbb{{\mathbb H}}
\def\Ubb{{\mathbb U}}
\def\Nbb{{\mathbb N}}

\def\zvec{\vec{z}}
\def\xvec{\vec{x}}
\def\zvecp{\vec{z}^{\prime}}


\title[Global estimates]%
{Global estimates for the Hartree-Fock-Bogoliubov equations}

\author{J. Chong}
\address{University of Texas at Austin}
\email{jwchong@math.utexas.edu}

\author{M. Grillakis}
\address{University of Maryland, College Park}
\email{mng@math.umd.edu}

\author{M. Machedon}
\address{University of Maryland, College Park}
\email{mxm@math.umd.edu}

\author{Z. Zhao}
\address{University of Maryland, College Park}
\email{zzh@umd.edu}

\date{\today}
\dedicatory{}
\commby{}

\maketitle
\begin{abstract}
We prove that certain Sobolev-type norms, slightly stronger than those given by energy conservation, stay bounded uniformly in time and $N$. This allows one to extend the local existence results of the second and third author globally in time. The proof is based on interaction Morawetz-type estimates and Strichartz estimates (including some new end-point results)
 for the equation
$\{ \frac{1}{i}\partial_t-\Delta_{x}-\Delta_{y}+\frac{1}{N}V_N(x-y) \}\Lambda(t, x, y) =F$
 in mixed coordinates such as
$L^p(dt) L^q(dx) L^2(dy)$,
$L^p(dt) L^q(dy) L^2(dx)$,\\
 $L^p(dt) L^q(d(x-y)) L^2(d(x+y))$
. The main new technical ingredient is a dispersive estimate in mixed coordinates, which may be of interest in its own right.
\end{abstract}

\section{Introduction}

\bigskip

 This paper is devoted to the study of some global estimates for solutions to a coupled system of Schr\"odinger-type equations
 (see \eqref{phieq}, \eqref{leq} and \eqref{geq} below)
 approximating the evolution of weakly interacting Bosons.
 For the sake of completeness, we include a brief overview
of the argument motivating these equations.

We refer to \cite{G-M2017} for  detailed explanations.
The problem is to understand the linear Schr\"odinger evolution of data equal to (or close to) a tensor product $\phi(x_1) \cdots \phi(x_N)$. The Hamiltonian is
\begin{align*}
H_{PDE}=\sum_{j=1}^{N}\Delta_{x_{j}} - \frac{1}{N}
\sum_{i<j}V_N\big(x_{j}-x_{i}\big)
\end{align*}
(mean-field negative Hamiltonian). For simplicity, assume $V$ satisfies the following conditions
\begin{align}
&\mbox{$V$ is spherically symmetric and }\label{Vhyp} \\
&V \ge 0, \, V \in C_0^{\infty}, \,  \frac{\partial V}{\partial r}( r ) \le 0. \notag \label{Vhyp}
\end{align}
The problem is easier to understand in the symmetric Fock space, with  Hamiltonian
\begin{align*}
&\H:=\int dx\left\{
a^{\ast}_{x}\Delta a_{x}\right\} -\frac{1}{2N}\int dxdy\left\{V_{N}(x-y)a^{\ast}_{x}a^{\ast}_{y}a_{x}a_{x}\right\}.
\end{align*}
We recall that
the Fock space Hamiltonian acts as a PDE Hamiltonian on the $n$th entry of Fock space
\begin{align*}
H_{n, \, \, PDE}=\sum_{j=1}^{n}\Delta_{x_{j}} - \frac{1}{N}
\sum_{i<j}V_N(x_{i}-x_{j})
\end{align*}
(see for instance \cite{G-M2017} for the definition of Fock space and the creation and annihilation operators $a^*$ and $a$).
The natural choice for initial conditions is
\begin{align*}
 e^{-\sqrt{N}\A(\phi_{0})}e^{-\B(k_{0})}\Omega
 \end{align*}
where $\Omega$ is the Fock space vacuum,
\begin{align*}
\A(\phi):=\int dx\left\{\bar{\phi}(x)a_{x}-\phi(x)a^{\ast}_{x}\right\}
\end{align*}
so that $ e^{-\sqrt{N}\A(\phi)}$ is the Weyl (unitary) operator. The coherent state
\begin{equation*}
e^{-\sqrt{N}\A(\phi)}\Omega =\left(\ldots\ c_{n}\prod_{j=1}^{n}\phi(x_{j})\ \ldots\right)
\quad {\rm with}\quad  c_{n}=\big(e^{-N \|\phi\|^2_{L^2}}N^{n}/n!\big)^{1/2}\ .
\end{equation*}
has tensor products in each entry, making it a natural choice for this problem.

We also recall
 \begin{align*}
&\B(k):=
\frac{1}{2}\int dxdy\left\{\bar{k}(x,y)a_{x}a_{y}
-k(x,y)a^{\ast}_{x}a^{\ast}_{y}\right\}\ .
\end{align*}
The unitary operator $e^{-\B(k)}$ is called the implementation of a Bogoliubov transformation in the  Physics literature, and the Segal-Shale-Weil or the metaplectic representation in the Math literature.
The state
$e^{-\B(k)}\Omega$ is called a squeezed state in the Physics literature. It provides second-order corrections to coherent states.

In the recent math literature, this set-up first appeared
in \cite{Rod-S}, followed by \cite{GMM1} where $e^{-\B(k)}$ is formally introduced.

Thus the problem is to find ``effective equations" for $\phi$, $k$ so that the exact evolution

\begin{align}
 \psi_{exact}=e^{i t \H}e^{-\sqrt{N}\A(\phi_{0})}e^{-\B(k_{0})}\Omega \label{exact}
 \end{align}
 is approximated, in the Fock space norm, by the approximate evolution
 \begin{align}
 \psi_{approx}=e^{i \chi(t)}e^{-\sqrt{N}\A(\phi(t))}e^{-\B(k(t))}\Omega \ . \label{approx}
\end{align}
See \eqref{oldest} below for one such existing estimate.

The equations for $\phi$, $k$  are easier to understand in terms of $\phi$ and the auxiliary functions $\Lambda$ and $\Gamma$.
See \cite{G-M2017}.

We refer to \cite{BCS} for a result of this type, in a slightly different setting.
That work is not based on the coupled equations \eqref{phieq}, \eqref{leq} and \eqref{geq}.

Fock space techniques can also be applied to $L^2(\mathbb R^N)$ approximations. See the recent paper \cite{BNNS} and the references therein. We also mention the related approach of \cite{BPPS} and \cite{aBPPS}.
The equations we will study are similar in spirit to the Hartree-Fock-Bogoliubov equations for Fermions. For Bosons, they were derived
in \cite{GM2}, \cite{G-M2017}, and, independently, in \cite{BBCFS}  and the recent paper \cite{BSS}. The first two references treat pure states, as described below,  while  last two  treat the case of mixed, quasi-free states. The PDEs are the same in both cases.
This ends our overview of the motivation, and we proceed  with the analysis of the equations.

The functions described by these PDEs are: the condensate $\phi(t, x)$  and the density matrices
\begin{align}
&\Gamma(t, x_1, x_2)=\frac{1}{N}\left(\shb \circ \sh\right)(t, x_1, x_2)+ \bar \phi(t, x_1) \phi(t, x_2)\label{def1Gamma}\\
 &\Lambda(t, x_1, x_2)=\frac{1}{2N}\sht (t, x_1, x_2) +\phi(t, x_1) \phi(t, x_2) \ . \label{def1Lambda}
 \end{align}
The pair excitation function $k$ is an auxiliary function, which does not explicitly appear in the system.

Let $V \in C_0^{\infty}(\mathbb R^3)$, $V \ge 0$, and denote $V_N(x-y)=N^{3\beta}V(N^{\beta}(x-y))$ be the potential, with $0\le \beta\le1$. We consider the following system:

\begin{align}
    &\{ \frac{1}{i}\partial_t-\Delta_{x_1} \}\phi(t, x_1)=-\int \phi(x_1)V_N(x_1-y)\Gamma(y,y) dy \label{phieq}\\
    &- \int \{V_N(x_1-y)\phi(y)(\Gamma(y,x_1)-\bar{\phi}(y)\phi(x_1))+V_N(x_1-y)\bar{\phi}(y)(\Lambda(x_1,y)-\phi(x_1)\phi(y))\} dy ,\notag
\end{align}

\begin{align}
    &\{ \frac{1}{i}\partial_t-\Delta_{x_1}-\Delta_{x_2}+\frac{1}{N}V_N(x_1-x_2) \}\Lambda(t, x_1,x_2)  \label{leq} \\
    &=-\int \{ V_N(x_1-y)\Gamma(y,y)+V_N(x_2-y)\Gamma(y,y) \} \Lambda(x_1,x_2) dy \notag\\
    &\quad - \int \{(V_N(x_1-y)+V_N(x_2-y))(\Lambda(x_1,y)\Gamma(y,x_2)+\bar{\Gamma}(x_1,y)\Lambda(y,x_2))\} dy \notag \\
    &\quad +2\int \{ (V_N(x_1-y)+V_N(x_2-y)) |\phi(y)|^2 \phi(x_1)\phi(x_2) \} dy ,\notag
\end{align}

\begin{align}
      &\{ \frac{1}{i}\partial_t-\Delta_{x_1}+\Delta_{x_2}\}\bar{\Gamma}(t, x_1,x_2)\label{geq}\\
    &=-\int \{ (V_N(x_1-y)-V_N(x_2-y))\Lambda(x_1,y)\bar{\Lambda}(y,x_2) \}  dy \notag\\
    &\quad - \int \{(V_N(x_1-y)-V_N(x_2-y))(\bar{\Gamma}(x_1,y)\bar{\Gamma}(y,x_2)+\bar{\Gamma}(y,y)\bar{\Gamma}(x_1,x_2))\} dy\notag \\
    &\quad +2\int \{ (V_N(x_1-y)-V_N(x_2-y)) |\phi(y)|^2 \phi(x_1)\bar{\phi}(x_2) \} dy.\notag
\end{align}

The solutions $\phi, \Lambda$, and $\Gamma$ also depend on $N$. This has been suppressed to simplify the notation. However, we will always keep track of dependence on $N$ in our estimates.

In order to motivate our main result (Theorem \ref{goal} below),  we recall the conserved quantities of the system, which will also be used in the proof of our main theorem.

The first conserved quantity is the total number of particles (normalized by division by $N$)
and it is
\begin{align}
{\rm tr}\left\{\Gamma(t)\right\}= \|\phi(t, \cdot)\|^2_{L^2(dx)} +\frac{1}{N} \|\sh (t, \cdot, \cdot)\|^2_{L^2(dxdy)}=1\ . \label{b5-nbr0}
\end{align}
From here we see that
\begin{align}
\|\Lambda(t, \cdot, \cdot) \|_{L^2(dx dy)} \le C\ . \label{L^2cons}
\end{align}

The second conserved quantity is the energy per  particle
\begin{align}
E(t):=&\ {\rm tr}\left\{\nabla_{x_{1}}\cdot\nabla_{x_{2}}\Gamma(t)\right\}
+\frac{1}{2}
\int dx_{1}dx_{2}\left\{V_{N}(x_{1}-x_{2})\big\vert\Lm(t,x_{1},x_{2})\big\vert^{2}
\right\}
\label{energy}
\\
&+\frac{1}{2}
\int dx_{1}dx_{2}\left\{
V_{N}(x_{1}-x_{2})\left(\big\vert\Gamma(t,x_{1},x_{2})\big\vert^{2}
+\Gamma(t,x_{1},x_{1})\Gamma(t,x_{2},x_{2})\right)
\right\}
\nonumber
\\
&-\int dx_{1}dx_{2}\left\{V_{N}(x_{1}-x_{2})
\vert\phi(t,x_{1})\vert^{2}\vert\phi(t,x_{2}\vert^{2}
\right\}\ .
\notag
\end{align}

Of special interest is the
kinetic part of the energy ,
\begin{align}
&{\rm tr}\left\{\nb_{x_{1}}\cdot\nb_{x_{2}}\Gamma\right\}
=\int dx\left\{\vert\nb_{x}\phi(t,x)\vert^{2}\right\} \label{energy1}\\
&+\frac{1}{2N}\int dx_{1}dx_{2}\left\{\vert\nb_{x_{1}}\sh(t, x_1, x_2)\vert^{2}+\vert
\nb_{x_{2}}\sh(t, x_1, x_2)\vert^{2}\right\}\ . \notag
\end{align}

If we assume $E \le C$, then we have an $H^1$  estimate for $\Lambda$, uniformly in time (and $N$):
\begin{align}
&\int dx_{1}dx_{2}\left\{
\big\vert\nb_{x_{1}}\Lm\big\vert^{2}+\big\vert\nb_{x_{2}}\Lm(t, x_1, x_2)\big\vert^{2}\right\} \label{lambdaest}
\le C\\
&\mbox{and also} \notag\\
&\frac{1}{N}\int dx_{1}dx_{2}
\big\vert\nb_{x_1, x_2}\sht(t, x_1, x_2)\big\vert^{2}\notag
\le C.
\end{align}

 Also,
$\Gamma $
satisfies the $H^2$ type estimate
\begin{align*}
\big\Vert\vert\nb_{x_{1}}\vert\vert\nb_{x_{2}}\vert\Gamma(t)\Vert_{L^{2}(dx_{1}dx_{2})}
\leq E\ .
\end{align*}
See \cite{GM2}, \cite{G-M2017}, as well \cite{BBCFS}  for these conserved quantities.

In addition, we have an interaction Morawetz-type estimate:
if the initial conditions  have energy $\le C$ then
\begin{align*}
\|\phi(t, x)\|^2_{L^4(dt dx )}+ \|\Gamma(t, x, x)\|_{L^2(dt dx )} \le C\ .
\end{align*}

Recalling \eqref{def1Lambda}, we see right away that \eqref{lambdaest} can be improved (in different ways) for the two summands of $\Lambda$:
\begin{align}
&\int dx_{1}dx_{2}\left\{
\big\vert\nb_{x_{1}}\frac{1}{2N}\sht\big\vert^{2}+\big\vert\nb_{x_{2}}\frac{1}{2N}\sht(t, x_1, x_2)\big\vert^{2}\right\} \label{kest}
\le \frac{C}{N}\\
&\mbox{(decay in $N$) and} \notag\\
&\int dx_{1}dx_{2}
\big\vert\nb_{x_1}\phi(t, x_1)\nb_{x_2}\phi(x_2) \big\vert^{2}\notag
\le C \notag \\
&\mbox{(extra differentiablility).}\notag
\end{align}

The goal of this paper is to prove the following improvement to \eqref{lambdaest}:
\begin{theorem} \label{goal} Let $\phi= \phi_N(t, x)$, $\Lambda=\Lambda_N(t, x, y)$ and $\Gamma=\Gamma_N(t, x, y)$
given by \eqref{def1Gamma}, \eqref{def1Lambda}
be solutions of \eqref{phieq}, \eqref{leq}, \eqref{geq} with smooth data (but not necessarily smooth uniformly in $N$), satisfying
\begin{align*}
&{\rm tr}\left\{\Gamma(0)\right\} \le C\\
&E(0) \le C \, \, \mbox{ (see \eqref{energy} for the definition of $E(t)$)}\\
&\||\nabla_x||\nabla_y| \Lambda(0, x, y)\|_{L^2} \le C N
\end{align*}

Let $V $ satisfy \eqref{Vhyp},  and denote $V_N(x-y)=N^{3\beta}V(N^{\beta}(x-y))$, with $0\le\beta\le1$.
Then there exists $\epsilon >0$ such that
\begin{align}
&\int
\big\vert|\nb_x|^{\frac{1}{2} + \epsilon}|\nb_y|^{\frac{1}{2} + \epsilon}\Lm(t, x, y)\big\vert^{2}dx dy \label{lambdaestimp}
\le C
\end{align}
uniformly in $t$ and $N$.
\end{theorem}

This is significant because
in \cite{G-M2019}
it was shown that, for $0< \beta <1$, under suitable assumptions on $V$, for every $\epsilon >0$,
there exists $T_0 >0$ depending only on
\begin{align*}
&\|<\nabla_x>^{\frac{1}{2}+\epsilon}<\nabla_y>^{\frac{1}{2}+\epsilon}\Lambda(0, \cdot)\|_{L^2}+
\|<\nabla_x>^{\frac{1}{2}+\epsilon}<\nabla_y>^{\frac{1}{2}+\epsilon}\Gamma(0, \cdot)\|_{L^2}\\
&+
\|<\nabla_x>^{\frac{1}{2}+\epsilon}\phi(0, \cdot)\|_{L^2}
\end{align*}
such that the system is well-posed (in a certain norm) on $[0, T_0]$, see Theorem 3.3 and Corollary 3.4 in \cite{G-M2019}.
Thus, estimate \eqref{lambdaestimp} extends the estimates of \cite{G-M2019} globally in time.

The results of  \cite{G-M2019}  together with an estimate of the form
\begin{align}
&\int dx dy
\big\vert|\nb_x|^{\frac{1}{2} + \epsilon}|\nb_y|^{\frac{1}{2} + \epsilon}\Lm(t, x, y)\big\vert^{2} \label{Jackyest}
\le C(t)
\end{align}
(which is similar to \eqref{lambdaestimp}, except that the bound is allowed to grow sub-linearly in time) were used
in \cite{jacky2}
to give a Fock space approximation  of the form
\begin{align}
&\|\psi_{exact}-\psi_{approx}\|_{\F}:=\|e^{i t \H}e^{-\sqrt{N}\label{oldest}
\A(\phi_{0})}e^{-\B(k(0))}\Omega-
e^{i \chi(t)}e^{-\sqrt{N}\A(\phi(t))}e^{-\B(k(t))}\Omega\|_{\F}\\
&\le \frac{Ce^{P(t)}}{N^{\frac{1-\beta}{2} }}\notag
\end{align}
for a polynomial $P(t)$, and $0<\beta <1$.
(See \eqref{exact}, \eqref{approx} for the definitions.)
 It is expected that the estimates of the current paper will lead to a better Fock space approximation.
This will be done in future work by the first and last author.

In addition, it is of general interest to know if Soblov norms higher than those given by energy conservation grow in time.
This was first accomplished for the non-linear Schr\"odinger equation in \cite{Bourgain98}.

The proof of \eqref{lambdaestimp} is immediate if we interpolate between \eqref{kest} and the following
\begin{theorem} \label{D^2}  Under the assumptions of Theorem \ref{goal},
there exists $p$ such that
\begin{align}
&\int
\big\vert|\nb_x||\nb_y|\Lm(t, x, y)\big\vert^{2}dx dy \label{Dlambdaestimp}
\le C N^p
\end{align}
uniformly in time.
\end{theorem}
\begin{remark} The power $p$ we obtained is not optimal. However, it should be noted that, even if
$\big\vert|\nb_x||\nb_y|\Lm(t, x, y)\big\vert^{2}\le C$ at $t=0$, an estimate of this form (uniform in $N$) is not expected to hold at later times because of singularities induced by the potential $V_N$.
\end{remark}
The rest of this paper is devoted to  the proof of Theorem \ref{D^2}.
We regard the equation for $\Lambda$ as a linear equation with non-local ``coefficients" given by $\Gamma$ and a forcing term involving $\phi$. For $\Gamma$ and   $\phi$, we will only use {\it a priori} estimates, given by  conserved quantities and an interaction Morawetz estimate.

In addition, the proof involves new Strichartz estimates in mixed coordinates.

To give an idea of the proof, differentiating \eqref{leq},
\begin{align}
&\{ \frac{1}{i}\partial_t-\Delta_{x}-\Delta_{y}+\frac{1}{N}V_N(x-y) \}\nb_{x}\nb_{y}\Lambda(t, x, y) \notag\\
&=-(V_N\ast \Gamma(t,x,x)+V_N\ast \Gamma(t,y,y))\cdot \nb_{x}\nb_{y}\Lambda(t,x,y)\label{mainmorawetz}\\\
&+2\nb_{x}\nb_{y}\left(V_N * |\phi|^2(t, x)\phi(t, x)\phi(t, y)\right) + \, \mbox{ other terms}.  \notag
\end{align}
For the main term \eqref{mainmorawetz}, we divide the time interval $[0, \infty)$ into finitely many intervals (independent of $N$) such that $\|\Gamma(t, x, x)\|_{L^2(dt dx)}$ is small, and the contributions of this term can be absorbed in the left hand side.
This uses an idea of Bourgain \cite{Bourgain98} and an interaction Morawetz argument.
Based on the above conserved quantities and the interaction Morawetz estimate, it is easy to prove
\begin{align*}
\|\nb_{x}\nb_{y}\left(V_N * |\phi|^2(t, x)\phi(t, x)\phi(t, y)\right)\|_{L^2(dt) L^{\frac{6}{5}}(dx) L^2(dy)} \le C N^{power}\ .
\end{align*}
In fact, we will show that all the other remaining terms on the right-hand side are in a dual Strichartz space, with norms possibly growing in $N$. In order to show that, we will first have to
  estimate $\Lambda$ and $\nabla \Lambda$ in various Strichartz norms.

Then we get the desired result, provided we can  prove Strichartz estimates (including some  end-points) for the equation
\begin{align*}
&\{ \frac{1}{i}\partial_t-\Delta_{x}-\Delta_{y}+\frac{1}{N}V_N(x-y) \}\Lambda(t, x, y) =F.
\end{align*}
Proving these Strichartz estimates is the main new technical accomplishment of our current paper.

\bigskip

Acknowledgement. J. Chong was supported by the NSF through the RTG grant DMS- RTG 1840314.

\section{Strichartz estimates}
From now on we use the notation $A \lesssim B$ to mean: there exists $C$, independent of $N$, such that $A \le CB$.
\subsection{Set-up}
Let $(p_1, q_1)$, $(p_2, q_2)$ be Strichartz admissible pairs  in 3 space dimensions ($\frac{2}{p_i}+ \frac{3}{q_i}=\frac{3}{2}$), with  $p_i \ge 2$, and let $p_i'$, $q_i'$ the dual exponents.

Recall  $\frac{1}{N}V_N(x)= N^{ 3 \beta -1} V(N^{\beta }x)$, $0 \le \beta \le 1$.
Since the results of this section may be of general interest, we point out the properties of $V$ that will be used
(which are weaker than \eqref{Vhyp}).

We only assume $V \in L^{\frac{3}{2}}$, thus
$\frac{1}{N}V_N \in L^{\frac{3}{2}}$ uniformly in $N\ge 1$ and
$V(x)$ is such that we already know the homogeneous Strichartz estimate
\begin{align}
\| e^{i t(\Delta_x -\frac{1}{N}V_N(x))} f \|_{L^{p_1}(dt) L^{q_1}(dx)} \lesssim \| f\|_{L^2(dx)}\label{3+1homest}
\end{align}
uniformly in $N$, as well as the double end-point $3+1$ Strichartz inhomogeneous estimate
\begin{align}
\|\int_0^t e^{i (t-s)(\Delta_x -\frac{1}{N}V_N(x))} F(s) ds\|_{L^{p_1}(dt) L^{q_1}(dx)} \le C \| F\|_{L^{p_2'}(dt) L^{q_2'}(dx)}\label{3+1inhomest}
\end{align}
  with bounds
 independent of $N$.

 These assumptions hold for $V$ satisfying \eqref{Vhyp}:
  If $\beta <1$, just $V \in L^{\frac{3}{2}}$ and
 $N$ large is sufficient. In that case,  $\|\frac{1}{N}V_N\|_{L^{\frac{3}{2}}}$ is small and an easy perturbation argument proves
 \eqref{3+1homest}, \eqref{3+1inhomest}.

 If $\beta =1$,  and $V \in C_0^{\infty}$, $V \ge 0$, the estimates \eqref{3+1homest}, \eqref{3+1inhomest} follow by scaling from the
 corresponding estimates for $N=1$. In turn, these follow by the Keel-Tao \cite{K-T} argument from the dispersive estimate
 \begin{align*}
 \| e^{i t(\Delta_x -V)} f \|_{L^{\infty}(\mathbb R^3)} \lesssim \frac{1}{t^{\frac{3}{2}}} \|f\|_{L^1(\mathbb R^3)}.
 \end{align*}
 There is an extensive literature on such estimates, following the breakthrough paper \cite{JSS}, but we could not find an explicit discussion of the case
 $V \in C_0^{\infty}(\mathbb R^3)$, $V \ge 0$. However, this follows, for instance, from \cite{yajima}, Theorem 1.3\footnote{
 In fact, just part 2 of Lemma 2.2 in \cite{yajima} suffices to prove the Strichartz estimates \eqref{3+1homest}, \eqref{3+1inhomest}, by standard Kato smoothing techniques. This avoids using the harder dispersive estimate.
 }. Since $-\Delta_x +V$ is a non-negative operator, it has no negative eigenvalues.  It is well-known
 $-\Delta_x +V$
  has no positive eigenvalues
(by Kato's theorem \cite{Kato}, or the earlier and more elementary result \cite{Rellich}, for instance).
 It is easy to show that $0$ is not a resonance or eigenvalue. The corresponding solution to  $(-\Delta_x +V)u=0$ is harmonic away from the support of $V$ and, if $u$ satisfies the resonance condition $<x>^{- \gamma} u \in L^2$ for all $\gamma > 1/2$, then, using the mean-value theorem one gets  $|u(x)| \lesssim |x|^{\gamma - \frac{3}{2}}$,  $|\nabla u(x)| \lesssim |x|^{\gamma - \frac{5}{2}}$ for $|x|$ sufficiently large. Thus one can integrate by parts and get
 \begin{align*}
 \int |\nabla u|^2 + V|u|^2 =0
 \end{align*}
 thus $u=0$ and Theorem 1.3 in \cite{yajima}  can be applied.

 \subsection{Statement of the Strichartz estimate}

 The main results of this section
 refer to the equation
 \begin{align}
&\left(\frac{1}{i}\frac{\partial}{\partial t}  - \Delta_x - \Delta_y + \frac{1}{N} V_N(\frac{x-y}{\sqrt 2}) \right)\Lambda =F \label{maineq}\\
&\Lambda(0, x, y)=\Lambda_0(x, y). \notag
\end{align}

The natural Strichartz norm for our system of Hartree-Fock-Bogoliubov type equations is
\begin{align*}
\|\Lambda\|_{\mathcal{S}^{p, q}}= \max \big\{ \|\Lambda\|_{L^p(dt) L^q(dx) L^2(dy)},
\|\Lambda\|_{L^p(dt) L^q(dy) L^2(dx)}, \|\Lambda\|_{L^p(dt) L^q(d(x-y)) L^2(d(x+y))} \}
\end{align*}
with the dual Strichartz norm
\begin{align*}
&\|F\|_{\mathcal{S}_{dual}^{p', q'}}\\
&= \min \big\{ \|F\|_{L^{p'}(dt) L^{q'}(dx) L^2(dy)},
 \|F\|_{L^{p'}(dt) L^{q'}(dy) L^2(dx)},  \|F\|_{L^{p'}(dt) L^{q'}(d(x-y)) L^2(d(x+y))} \}
\end{align*}
and the natural question to ask is whether
\begin{align}
\|\Lambda\|_{\mathcal{S}^{p_1, q_1}} \lesssim  \|\Lambda_0\|_{L^2} + \|F\|_{\mathcal{S}_{dual}^{p_2', q_2'}}\label{slickStrichartz}
\end{align}
for any admissible  pairs $(p_1, q_1)$, $(p_2, q_2)$. This amounts to 9 inequalities. We will show that if not both $(p_1, q_1)$, $(p_2, q_2)$ are end-point exponents ($p=2, q=6$), then \eqref{slickStrichartz} is true (all 9 cases hold). In the double end-point case
we have to exclude the two cases where $x$ and $y$ are flipped:
we don't know if
\begin{align}
\|\Lambda\|_{L^2(dt) L^6(dx) L^2(dy)} \lesssim  \|\Lambda_0\|_{L^2}+ \|F\|_{L^2(dt) L^{6/5}(dy) L^2(dx)} \label{rotated??}
\end{align}
is true.

In order to exclude this, we fix a number $p_0>2$ (in our application,  $p_0=\frac{8}{3}, q_0=4$ will suffice) and define the "restricted" Strichartz norm
\begin{align}
&\|\Lambda\|_{\mathcal{S}_{restricted}}\label{restr}\\
&=  \sup_{p_0\le p \le \infty,\, \,  p, q \, \, admissible}\|\Lambda\|_{L^{p}(dt) L^{q}(dx) L^2(dy)}\notag\\
&\quad +\sup_{p_0\le p \le \infty, \, \,  p, q \, \notag \, admissible}\|\Lambda\|_{L^{p}(dt) L^{q}(dy) L^2(dx)} \\
 &\quad +
\sup_{2 \le p \le \infty, \, \,  p, q \, \, admissible}\|\Lambda\|_{L^p(dt) L^q(d(x-y)) L^2(d(x+y))}. \notag
\end{align}
Notice that the end-point is included in $x-y, \, x+y $ coordinates.

In this section, we prove
 \begin{theorem} (non-endpoint result) \label{nonendpointthm} Let $V\in L^{3/2}(\mathbb R^3)$ as above,  $0 \le \beta \le 1$
 and assume  \eqref{3+1homest},  \eqref{3+1inhomest} hold.
 Let $p_i, q_i$ ($i=1, 2$) be  Strichartz admissible pairs  and assume both $p_i >2$. Let $p_i'$, $q_i'$ be the dual exponents.
If $\Lambda$ satisfies \eqref{maineq}, then
\begin{align}
\|\Lambda\|_{\mathcal{S}^{p_1, q_1}} \lesssim  \|\Lambda_0\|_{L^2} + \|F\|_{\mathcal{S}_{dual}^{p_2', q_2'}}.\label{slickStrichartznopot}
\end{align}

\end{theorem}
We also have a ``one end-point result":
 \begin{theorem} (one endpoint result) \label{oneendpointthm} Let $V\in L^{3/2}$, $0 \le \beta \le 1$,
 and assume  \eqref{3+1homest},  \eqref{3+1inhomest} hold. Let $p_1, q_1$  be  Strichartz admissible pair  and assume  $p_1 >2$  or $p_2 >2$.
If $\Lambda$ satisfies \eqref{maineq}
then
\begin{align}
\|\Lambda\|_{\mathcal{S}^{p_1, q_1}} \lesssim  \|\Lambda_0\|_{L^2} + \|F\|_{\mathcal{S}_{dual}^{p_2', q_2'}}.\label{oneendpoint1}
\end{align}
\end{theorem}

Finally, we have a double end-point result:
\begin{theorem} \label{mainthm} Let $V\in L^{3/2}$, $0 \le \beta \le 1$ and  assume  \eqref{3+1homest},  \eqref{3+1inhomest} hold.

If $\Lambda$ satisfies \eqref{maineq},
then
\begin{align}
\|\Lambda\|_{\mathcal{S}^{2, 6}}
 \lesssim \|\Lambda_0\|_{L^2}+
\|F\|_{L^2(dt) L^{6/5}(d(x-y)) L^2(d(x+y))}.
 \label{unrotated}
\end{align}
\end{theorem}

\begin{remark}
The proof of the above theorem could be adapted to show the additional estimates
\begin{align*}
&\|\Lambda\|_{L^2(dt) L^6(dx) L^2(dy)} \lesssim  \|\Lambda_0\|_{L^2}+ \|F\|_{L^2(dt) L^{6/5}(dx) L^2(dy)} \\
&\|\Lambda\|_{L^2(dt) L^6(d(x-y)) L^2(d(x+y))} \lesssim  \|F\|_{\mathcal{S}_{dual}^{6/5, 2}}
\end{align*}
but, in order to keep the exposition simple, we won't do it.
\end{remark}
Theorem  \ref{oneendpointthm} and Theorem \ref{mainthm} imply the following concise form, which is what we will use in our  applications:

\begin{theorem}\label{concise} Let $V$ as above, $0 \le \beta \le 1$,  and $p_0>2$ defining $\mathcal{S}_{restricted}$ (see \eqref{restr}) be fixed.
If $\Lambda$ satisfies \eqref{maineq},
 then, for any admissible
   Strichartz   pair $(p, q)$ (including the end-point $(2, 6)$),
\begin{align}
\|\Lambda\|_{\mathcal{S}_{restricted}}
 \lesssim
 \|\Lambda_0\|_{L^2} + \|F\|_{\mathcal{S}_{dual}^{p', q'}}.
 \label{unrotated}
\end{align}
\end{theorem}

\begin{remark} The above theorems have immediate and obvious generalizations to all dimensions $\ge 3$. Also, the spaces can be localized to any finite or infinite time interval, and the theorems go through with obvious modifications. For instance,
\begin{align*}
\|\Lambda\|_{\mathcal{S}_{restricted}[T_1, T_2]}
 \lesssim
 \|\Lambda(T_1)\|_{L^2} + \|F\|_{\mathcal{S}_{dual}^{p', q'}[T_1, T_2]}.
\end{align*}
\end{remark}

\begin{remark} Obviously, Theorem \ref{oneendpointthm} implies Theorem \ref{nonendpointthm}. We list them separately because the proof of
Theorem \ref{nonendpointthm} is  based on standard techniques, while the proof of  Theorem \ref{oneendpointthm} and Theorem \ref{mainthm} requires essentially new ideas.
\end{remark}

 These are presented in the next two subsections.
\subsection{Standard techniques}
We will use the following well-known identities, which were also used  in \cite{Dancona}, \cite{Bouclet}, \cite{Mizutani}.

\begin{proposition} \label{Ndef}
Let
\begin{align*}
&\N F= i\int_0^t e^{i (t-s)(\Delta_x + \Delta_y -\frac{1}{N}V_N(\frac{x-y}{\sqrt 2}))} F(s) ds\\
&\N_0 F=i \int_0^t e^{i (t-s)(\Delta_x + \Delta_y )} F(s) ds.
\end{align*}
Then the following identities hold (denoting $V_N=V_N(\frac{x-y}{\sqrt 2})$)

\begin{align}
\N - \N_0 = - \N \frac{1}{N}V_N \N_0=- \N_0 \frac{1}{N}V_N \N\ .  \label{pot1}
\end{align}
and thus
\begin{align}
\N  =\N_0 - \N_0 \frac{1}{N}V_N \N_0 + \N_0 \frac{1}{N}V_N \N \frac{1}{N}V_N \N_0 . \label{pot2}
\end{align}
\end{proposition}

\begin{proof}
Look at
\begin{align*}
&\N \frac{1}{N}V_N \N_0=
\N \left(\left( \frac{1}{i}\frac{\partial}{\partial t}  - \Delta_x - \Delta_y + \frac{1}{N}V_N\right)-
\left( \frac{1}{i}\frac{\partial}{\partial t}  - \Delta_x - \Delta_y \right)\right) \N_0\\
&= \N_0-\N
\end{align*}
where we have used the fact that $\N$ and $\N_0$ are left and right inverses of the corresponding differential operators.
For the second part of \eqref{pot1}, reverse  the order of $\N$ and $\N_0$. The formula \eqref{pot2} is obtained by iterating
\eqref{pot1}.
\end{proof}

In addition, we need the following propositions:
\begin{proposition}\label{old1}
Let $\N_0$ be as in Proposition \ref{Ndef}. Let $(p_1, q_1)$, $(p_2, q_2)$ be Strichartz admissible (including the end-points $p_i=2, q_i=6$). Then
\begin{align}
&\|\N_0 F\|_{L^{p_1}(dt) L^{q_1}(dx) L^2(dy)} \lesssim \|F\|_{L^{p'_2}(dt) L^{q_2'}(dx) L^2(dy)}\label{est1}\\
&\big\|e^{i t (\Delta_x + \Delta_y)} \Lambda_0\big\|_{L^{p_1}(dt) L^{q_1}(dx) L^2(dy)} \lesssim \|\Lambda_0\|_{L^2}.\label{est2}
\end{align}
\end{proposition}
\begin{proof}

 \begin{align*}
 \|\N_0 F(t, x, \cdot)\|_{L^2(dy)} = &\ \|e^{i t \Delta_y} \int_0^t e^{i(t-s)  \Delta_x}e^{- i s  \Delta_y}F (s, \cdot, \cdot) ds\|_{L^2(dy)}\\
 =&\ \| \int_0^t e^{i(t-s)  \Delta_x}e^{- i s  \Delta_y}F (s, \cdot, \cdot) ds\|_{L^2(dy)}\\
 \end{align*}
 and
 \begin{align*}
\big\| \|\N_0 F(t, x, \cdot)\|_{L^2(dy)}\big \|_{L^{p_1}(dt)L^{q_1}(dx)} &= \big\| \int_0^t e^{i(t-s)  \Delta_x}e^{- i s  \Delta_y}F (s, \cdot, \cdot) ds\|_{L^2(dy)}\big\|_{L^{p_1}(dt)L^{q_1}(dx)}\\
& \le \big\|\| \int_0^t e^{i(t-s)  \Delta_x}e^{- i s  \Delta_y}F (s, \cdot, \cdot) ds\|_{L^{p_1}(dt)L^{q_1}(dx)}\big\| _{L^2(dy)}\\
&\le C \big\| \|e^{- i s  \Delta_y}F (s, \cdot, \cdot)\|_{L^{p'_2}(dt) L^{q_2'}(dx)}\big\| _{L^2(dy)}\\
&\le C \big\| \|e^{- i s  \Delta_y}F (s, \cdot, \cdot)\|_{L^2(dy)}\big\| _{L^{p'_2}(dt) L^{q_2'}(dx)}\\
&=C \|F\|_{L^2(dt) L^{6/5}(dx) L^2(dy)}.
\end{align*}
The proof of \eqref{est2} is similar. See Lemma 5.3 in \cite{G-M2017}.
\end{proof}

We also have the following version which excludes the double end-point, but works with any choice of coordinate systems:

\begin{proposition}\label{Ch-K}
Let $\N_0$ be as in Proposition \ref{Ndef}. Let $p_i, q_i$ ($i=1, 2$) be  Strichartz admissible pairs, with at least one $p_i>2$. Also, let $R \in O(6)$.  Then
\begin{align*}
&\|\N_0 F\|_{L^{p_1}(dt) L^{q_1}(dx) L^2(dy)} \lesssim \|F \circ R\|_{L^{p_2'}(dt) L^{q_2'}(dx) L^2(dy)}.
\end{align*}
In particular,
\begin{align}
\|N_0 F\|_{\mathcal{S}^{p_1, q_1}} \lesssim   \|F\|_{\mathcal{S}_{dual}^{p_2', q_2'}}.\label{Christ}
\end{align}
\end{proposition}

\begin{proof} Using \eqref{est2}, the $T T^*$ argument and the $O(6) $ invariance of $\Delta$ we have
\begin{align*}
\big\|\int_0^\infty e^{i(t-s)(\Delta_x + \Delta_y)} F(s, \cdot) ds\big\|_{L^{p_1}(dt) L^{q_1}(dx) L^2(dy)} \lesssim \|F \circ R\|_{L^{p_2'}(dt) L^{q_2'}(dx) L^2(dy)}.
\end{align*}
By the Christ-Kiselev lemma (Lemma 2.4 in \cite{Tao}),  we conclude
\begin{align*}
\big\|\int_0^t e^{i(t-s)(\Delta_x + \Delta_y)} F(s, \cdot) ds\big\|_{L^{p_1}(dt) L^{q_1}(dx) L^2(dy)} \lesssim \|F \circ R\|_{L^{p_2'}(dt) L^{q_2'}(dx) L^2(dy)}
\end{align*}
provided $p_1>p_2'$.
\end{proof}
Finally, we have a version which includes the potential, but only works in coordinates
compatible with the potential:
\begin{proposition} \label{samesplitting}
 If $V(x)$ is such that we already know  \eqref{3+1homest},  \eqref{3+1inhomest}.
Then,
\begin{align}
&\|e^{i t \left( \Delta_x + \Delta_y -\frac{1}{N}V_N(\frac{x-y}{\sqrt 2})\right)}\Lambda_0\|_{L^{p_1}(dt) L^{q_1}(d(x-y)) L^2(d(x+y))} \lesssim \|\Lambda_0\|_{L^2(dx dy)}\label{est3}\\
&\|\N F\|_{L^{p_1}(dt) L^{q_1}(d(x-y)) L^2(d(x+y))} \lesssim \|F\|_{L^{p_2'}(dt) L^{q_2'}(d(x-y)) L^2(d(x+y))}. \label{double1}
\end{align}
\end{proposition}

\begin{proof} The proof is similar to that of \eqref{est1} and \eqref{est2}, but is based on writing $  \Delta_x + \Delta_y -V_N(x-y)=
 \Delta_{\frac{x+y}{\sqrt 2}} +\left( \Delta_{\frac{x-y}{\sqrt 2}} - V_N(\frac{x-y}{\sqrt 2})\right)
$ and using the fact that these commute.
\end{proof}

\subsection{The new estimate}

The main step in the end-point cases,  which may be of interest in its own right,  does not involve the potential. We will show

\begin{theorem} \label{new}
Let $\Lambda= \N_0 F$ be the solution to
\begin{align*}
&\left(\frac{1}{i}\frac{\partial}{\partial t}  - \Delta_x - \Delta_y  \right)\Lambda =F\\
&\Lambda(0, x, y)=0.
\end{align*}
Then the following closely related estimates hold:
\begin{align}
&\|\Lambda\|_{L^2(dt) L^6(dx) L^2(dy)} \le C \|F\|_{L^2(dt) L^{6/5}(d(x-y)) L^2(d(x+y))} \label{mix}\\
&\|\Lambda\|_{L^2(dt) L^6(dy) L^2(dx)} \le C \|F\|_{L^2(dt) L^{6/5}(d(x-y)) L^2(d(x+y))} \label{mix2}\\
&\mbox{and also,}\notag\\
&\|\Lambda\|_{L^2(dt) L^6(d(x-y)) L^2(d(x+y))} \le C \|F\|_{L^2(dt) L^{6/5}(dx)) L^2(dy)}\label{mix1}\\
&\|\Lambda\|_{L^2(dt) L^6(d(x-y)) L^2(d(x+y))} \le C \|F\|_{L^2(dt) L^{6/5}(dy) L^2(dx)}. \notag
\end{align}

\end{theorem}

Together with the estimates of the previous subsection, Theorem \ref{new} implies
\begin{corollary} For any Strichartz admissible pair $p, q$ (including the end-point)
\begin{align}
&\|N_0 F\|_{\mathcal{S}^{p, q}} \lesssim   \|F\|_{L^2(dt) L^{6/5}(d(x-y)) L^2(d(x+y))}\\
&\|N_0 F\|_{L^2(dt) L^6(d(x-y)) L^2(d(x+y))} \lesssim  \|F\|_{\mathcal{S}_{dual}^{p', q'}}.
\end{align}
\end{corollary}
This complements the estimates of  Proposition \ref{old1}, Proposition \ref{samesplitting}, and Proposition \ref{Ch-K}. And, it will be used in the proof
of Theorem \ref{mainthm}.

The proof of Theorem \ref{new} will be given in subsection \ref{newproof}. It uses
 a new dispersive estimate in mixed coordinates, see
Proposition \ref{dispersive} below.

Now we can outline the proofs of our main results.
\subsection{Proofs of Theorem \ref{nonendpointthm},  Theorem \ref{oneendpointthm} and
Theorem \ref{mainthm}, assuming Theorem \ref{new}.}

\begin{proof}
   Assume first $\Lambda_0=0$.
We proceed to estimate the terms in \eqref{pot2}.

\begin{align*}
\N F =\N_0 F - \N_0 \frac{1}{N}V_N \N_0 F + \N_0  \frac{1}{N}V_N \N  \frac{1}{N}V_N \N_0 F.
\end{align*}

For the first term, if $p_1 >2$ or $p_2 >2 $
use Proposition \ref{Ch-K}:
\begin{align*}
&\|\N_0 F\|_{\mathcal S{p_1, q_1}} \lesssim
\|F\|_{\mathcal{S}_{dual}^{p_2', q_2'}}
\end{align*}
while, for the proof of Theorem \ref{mainthm}, if we are in the double end-point case, we use Theorem  \ref{new}:
\begin{align*}
\|\N_0 F\|_{\mathcal{S}^{2, 6}}
 \lesssim
\|F\|_{L^2(dt) L^{6/5}(d(x-y)) L^2(d(x+y))}.
\end{align*}

 This is the only term where we don't know if we can flip $x$ and $y$ in the double end-point case.

For the second term,
\begin{align*}
&\|\N_0  \frac{1}{N}V_N \N_0 F\|_{\mathcal S^{p_1, q_1}} \lesssim
\| \frac{1}{N}V_N \N_0 F\|_{L^{2}(dt) L^{\frac{6}{5}}(d(x-y)) L^2(d(x+y))}\\ &(\mbox{we used Proposition \ref{Ch-K} if  $p_1>2$ and Theorem  \ref{new} if $p_1=2$})\\
&\lesssim \| \frac{1}{N}V_N\|_{L^{\frac{3}{2}}}\|\N_0 F\|_{L^{2}(dt) L^{6}(d(x-y)) L^2(d(x+y))}.
\end{align*}
Using Proposition \ref{Ch-K} if $p_2>2$
and Theorem \ref{new} if $p_2=2$,
 we conclude
\begin{align*}
& \|\N_0 F\|_{L^{2}(dt) L^{6}(d(x-y)) L^2(d(x+y))}\lesssim \|F\|_{\mathcal{S}_{dual}^{p_2', q_2'}}.
\end{align*}

For the third term in \eqref{pot2} we proceed along the same lines,
\begin{align*}
&\|\N_0  \frac{1}{N}V_N \N  \frac{1}{N}V_N \N_0 F\|
_{\mathcal S^{p_1, q_1}} \\
& \lesssim
\| \frac{1}{N}V_N \N  \frac{1}{N}V_N \N_0 F\|_{L^{2}(dt) L^{\frac{6}{5}}(d(x-y)) L^2(d(x+y))}\\
&\lesssim \|\N  \frac{1}{N}V_N \N_0 F\|_{L^{2}(dt) L^6(d(x-y)) L^2(d(x+y))}\\
&\lesssim
\|\frac{1}{N} V_N \N_0 F\|_{L^{2}(dt) L^{\frac{6}{5}}(d(x-y)) L^2(d(x+y))}\\
&\mbox{(here we used Proposition \ref{samesplitting})}\\
&\lesssim \| \frac{1}{N}V_N\|_{L^{\frac{3}{2}}}\|\N_0 F\|_{L^{2}(dt) L^{6}(d(x-y)) L^2(d(x+y))}\\
&\lesssim \|F\|_{\mathcal{S}_{dual}^{p_2', q_2'}}.
\end{align*}
Notice that if either $p_1=2$ or $p_2=2$ we have to use Theorem \ref{new}.

Finally, we show how to reduce the proof of  Theorem \ref{nonendpointthm},  Theorem \ref{oneendpointthm} and
Theorem \ref{mainthm} to the case $\Lambda_0=0$. Consider the homogeneous version of the above Theorems ($F=0$),
written in the form
\begin{align*}
&\left(\frac{1}{i}\frac{\partial}{\partial t}  - \Delta_x - \Delta_y  \right)\Lambda =-  \frac{1}{N}V_N(\frac{x-y}{\sqrt 2}) \Lambda\\
&\Lambda(0, x, y)=\Lambda_0,
\end{align*}
where we treat $ \frac{1}{N}V_N(\frac{x-y}{\sqrt 2}) \Lambda= \frac{1}{N}V_N(\frac{x-y}{\sqrt 2})e^{i t (\Delta_x + \Delta_y- \frac{1}{N}V_N(\frac{x-y}{\sqrt 2}))} \Lambda_0 $ as a forcing term.

 From \eqref{est3}  we have, for
 $\Lambda= e^{i t (\Delta_x + \Delta_y- \frac{1}{N} V_N(\frac{x-y}{\sqrt 2}))} \Lambda_0$,
\begin{align*}
&\big\|e^{i t (\Delta_x + \Delta_y-  \frac{1}{N}V_N(\frac{x-y}{\sqrt 2}))} \Lambda_0\big\|_{L^2(dt) L^6(d(x-y)) L^2(d(x+y))} \lesssim \|\Lambda_0\|_{L^2}
\end{align*}
thus
\begin{align*}
&\| \frac{1}{N}V_N(\frac{x-y}{\sqrt 2})\Lambda\|_{L^2(dt)L^{6/5}(d(x-y))L^2(d(x+y))}\\
&\lesssim \| \frac{1}{N}V_N\|_{L^{3/2}}\|\Lambda\|_{L^2(dt)L^6(d(x-y))L^2(d(x+y))} \lesssim  \|\Lambda_0\|_{L^2}
\end{align*}
and we use Proposition \ref{Ch-K} or
Theorem \ref{new}
to conclude
\begin{align*}
&\|\N_0 \left( \frac{1}{N}V_N(\frac{x-y}{\sqrt 2}) \Lambda\right)\|_{\mathcal S^{p_1, q_1}}\\
&\lesssim \| \frac{1}{N}V_N(\frac{x-y}{\sqrt 2})\Lambda\|_{L^2(dt)L^{6/5}(d(x-y))L^2(d(x+y))}
 \lesssim  \|\Lambda_0\|_{L^2}.
\end{align*}

Finally, from \eqref{est2} we have
\begin{align*}
&\big\|e^{i t (\Delta_x + \Delta_y)} \Lambda_0\big\|_{L^2(dt) L^6(dx) L^2(dy)} \lesssim \|\Lambda_0\|_{L^2}.
\end{align*}

\end{proof}

It remains to prove Theorem \ref{new}.
\subsection{Proof of Theorem \ref{new}}\label{newproof}

The proof will follow the outline of Keel and Tao.
The main step is  proving a new dispersive estimate.

\begin{proposition} \label{dispersive}
\begin{align*}
\|e^{it \left(\Delta_x + \Delta_y\right)} f\|_{L^{\infty}(d(x-y)) L^2(d(x+y))} \le \frac{C}{t^{3/2}}  \|f\|_{L^1(dx) L^2(dy)}
\end{align*}
and, similarly,
\begin{align}
\|e^{it \left(\Delta_x + \Delta_y\right)} f\|_{L^{\infty}(dx) L^2(dy)}
\label{inv}
 \le \frac{C}{t^{3/2}}  \|f\|_{L^1(d(x-y)) L^2(d(x+y))}.
\end{align}
\end{proposition}

\begin{proof}

Our proof is inspired, in part, by Lemma 1 in \cite{FLLS} and also Lemma 2.2 in \cite{JSS}.

We will prove \eqref{inv}.

By a density argument, it suffices to take
take
\begin{align*}
 f (x, y)= \sum u_k\left(\frac{y-x}{\sqrt 2}\right) v_k\left(\frac{x+y}{\sqrt 2}\right)
  \end{align*}
  with $u_k$ orthogonal (but not normalized), and $v_k$ orthonormal. (This is a singular value decomposition of $f$ composed with a rotation; it will turn out that the orthogonality of $u_k$ will not play a role).
 Then
 \begin{align*}
e^{it \left(\Delta_x + \Delta_y\right)} f(x, y) =
\sum \left(e^{it \Delta}u_k\right)\left(\frac{y-x}{\sqrt 2}\right)\left(e^{it \Delta} v_k\right)\left(\frac{x+y}{\sqrt 2}\right).
\end{align*}

Then the LHS of \eqref{inv} is $\sup_{x_0}\|\sum \left(e^{it \Delta}u_k\right)\left( \frac{\cdot-x_0}{\sqrt 2}\right)\left(e^{it \Delta}v_k\right) \left(\frac{\cdot+x_0}{\sqrt 2}\right)\|_{L^2(\mathbb R^3)}$. Look at this expression with $x_0$ fixed.

The RHS of \eqref{inv} is, using Plancherel and the fact that $v_k$ are orthonormal,
$RHS$ of  \eqref{inv}$ =\frac{C}{t^{3/2}} \|\left( \sum |u_k|^2\right)^{\frac{1}{2}}\|_{L^1(\mathbb R^3)}$.
The proof will be complete once we prove the following lemma, in which the general orthonormal set $\left(e^{it \Delta}v_k\right) (\frac{\cdot+x_0}{\sqrt 2})$ is re-labeled $v_k$ and the $u_k$ have also been shifted by $x_0$ and re-scaled by $\frac{1}{\sqrt 2}$ .
 \end{proof}

\begin{lemma} There exists $C>0$ such that,
for any $u_k$,
\begin{align}
\sup_{v_k \,  \mbox{ orthonormal}}\big\|\sum \left( e^{i t \Delta} u_k\right) v_k\big\|_{L^2(\mathbb R^3)} \le \frac{C}{t^{3/2}} \big\|\left(\sum|u_k|^2\right)^{\frac{1}{2}}\big\|_{L^1(\mathbb R^3)}. \label{gen}
\end{align}
\end{lemma}
\begin{proof}

Since we take supremum over all orthonormal sets  $v_k$, and $t$ is fixed, we may replace $v_k$ by $e^{-i t \Delta} v_k$, and \eqref{gen} is equivalent to
\begin{align}
\sup_{v_k \,  \mbox{ orthonormal}}\big\|\sum \overline{e^{i t \Delta} u_k}(x)e^{i t \Delta} v_k(x)\|_{L^2} \le \frac{C}{t^{3/2}} \big\|\left(\sum|u_k|^2\right)^{\frac{1}{2}}\big\|_{L^1}. \label{conj gen}
\end{align}

For any $A \in \mathcal S(\mathbb R^3)$,  let
$e^{-i t \Delta} A(x) e^{i t \Delta}=
A(x+2tD)$ where $D =p= \frac{1}{i} \frac{\partial}{\partial x}$.
Using the well-known formula
$$e^{-i t \Delta}e^{i x \cdot \xi} e^{i t \Delta}f(x) = e^{i x \cdot \xi} e^{i t |\xi|^2}f(x + 2 t \xi) $$
we compute
\begin{align*}
&e^{-i t \Delta} A(x) e^{i t \Delta} f(x) = \frac{1}{(2 \pi)^3} \int \hat A (\xi)e^{-i t \Delta}e^{i x \cdot \xi} e^{i t \Delta}f (x) d \xi\\
&=\frac{1}{(2 \pi)^3} \int \hat A(\xi) e^{i \xi \cdot x} e^{i t |\xi|^2}f(x + 2 t \xi) d \xi\\
&(\mbox{ change variables }\, \xi \to \frac{\xi-x}{2t})\\
&=\frac{1}{(4  \pi t)^3} \int \hat A\left(\frac{\xi-x}{2t}\right) e^{i \frac{\xi-x}{2t} \cdot x} e^{i t|\frac{\xi-x}{2t}|^2} f (\xi) d \xi\\
&=\frac{1}{(4  \pi t)^3} \int\hat A\left(\frac{\xi-x}{2t}\right) e^{-i \frac{|x|^2}{4t}} e^{i \frac{|\xi|^2}{4t}} f (\xi) d \xi.
\end{align*}

Thus  the integral kernel corresponding to
$A(x+ 2tD)$ is
\begin{align*}
 &K_t(x, y)=\frac{1}{(4  \pi t)^3} \hat A\left(\frac{-x+ y}{2t}\right) e^{-i \frac{|x|^2}{4t}} e^{i \frac{|y|^2}{4t}}\\
 &=B_{t, x}(y)  e^{-i \frac{|x|^2}{4t}} e^{i \frac{|y|^2}{4t}}
 \end{align*}
 where, in order to simplify the notation, for fixed $t, x$, we defined $B_{t, x}(y)=\frac{1}{(4  \pi t)^3}\hat A\left(\frac{-x+ y}{2t}\right)$.
Notice
\begin{align*}
\|B_{t, x}\|_{L^2(dy)}= \frac{c}{t^{\frac{3}{2}}} \|A\|_{L^2}.
\end{align*}

For a suitable $A$ with $\|A\|_{L^2}=1$,
\begin{align}
&\big\|\sum \overline{e^{i t \Delta} u_k}(x)e^{i t \Delta} v_k(x)\|_{L^2}\\
&=\int \sum \overline{e^{i t \Delta} u_k}(x) A(x)e^{i t \Delta} v_k(x) dx\notag\\
&=\sum <e^{i t \Delta}u_k, A e^{i t \Delta}v_k >=\sum <u_k,e^{-i t \Delta} A e^{i t \Delta}v_k >\notag\\
&= \sum <u_k, A(x +2 t D)v_k >.\label{continue}
\end{align}

From now we take any  $A \in \mathcal S(\mathbb R^3)$ with $\|A\|_{L^2(\mathbb R^3)}=1$.

We have to show
\begin{align*}
&|\eqref{continue}|=\bigg| \sum \int \overline{ e^{i\frac{|x|^2}{4t}}u_k(x)} B_{t, x}(y) e^{i \frac{|y|^2}{4t}}v_k(y)dx \, dy\bigg| \le \frac{C}{t^{\frac{3}{2}}}
\big\|\left(\sum|u_k|^2\right)^{\frac{1}{2}}\big\|_{L^1}
\end{align*}
for any orthonormal $v_k$ and any $\|A\|_{L^2(\mathbb R^3)}=1$.
The exponentials play no role now (change notation and remove them).

Look at

\begin{align*}
& \sum \int\overline{ u_k(x)}\left(\int  B_{t, x}(y)v_k(y)dy\right) \, dx\\
&= \sum \int \overline{u_k(x)}c_k(t, x) \, dx
\end{align*}
where, for fixed $t$ and $x$,
\begin{align*}
c_k(t, x)=\int   B_{t, x}(y)v_k(y)dy
\end{align*}
is a Fourier coefficient of $ B_{t, x}$. By Plancherel, we have $\sum |c_k(t, x)|^2 \le \|B_{t, x}\|_{L^2}^2$ uniformly in $t, x$.

Now we go back to
\begin{align*}
&\big|\sum \int \overline{ u_k(x)}c_k(t, x) \, dx\big|
\le \int \left( \sum |u_k(x)|^2\right)^{\frac{1}{2}}\left(\sum |c_k(t, x)|^2\right)^{\frac{1}{2}} dx\\
& \le\|B_{t, x}\|_{L^2(dy)}\big\|\left(\sum|u_k|^2\right)^{\frac{1}{2}}\big\|_{L^1}= \frac{c}{t^{\frac{3}{2}}} \|A\|_{L^2}
\big\|\left(\sum|u_k|^2\right)^{\frac{1}{2}}\big\|_{L^1}.
\end{align*}
A second proof of this proposition will be given in section \ref{secondproof}.

\end{proof}

We will finish the proof of Theorem \ref{new} by adapting the argument of Keel and Tao, \cite{K-T}.

Let  $R$ be the rotation
 $(x, y) \to \frac{1}{\sqrt 2} (x-y, x+y)$.
Following \cite{K-T},  define
\begin{align*}
T(F, G)=\int_{- \infty}^{\infty} \int_0^t <e^{i (t-s)\Delta_{x, y}} F(s), G\circ R (t)>ds dt
\end{align*}
with $T_j$ the above integral restricted to $t-2^{j+1}<s< t-2^j$.
In this formulation, the goal is $|T(F, G)| \le C \|F\|_{L^2(dt) L^{\frac{6}{5}}(dx) L^2 (dy)}  \|G\|_{L^2(dt) L^{\frac{6}{5}}(dx) L^2 (dy)}$.

Using the dispersive estimate of Proposition \ref{dispersive},  Lemma 4.1 in \cite{K-T} goes through word by word, and we have
\begin{align*}
|T_j(F, G)| \le C 2^{- j \beta(a, b)}\|F\|_{L^2(dt) L^{a'}(dx) L^2 (dy)}\|G\|_{L^2(dt) L^{b'}(dx) L^2 (dy)}
\end{align*}
for all $\left(\frac{1}{a}, \frac{1}{b}\right)$ in a neighborhood of $\left(\frac{1}{6}, \frac{1}{6}\right)$. Here $\beta(a, b)
=\frac{1}{2}-\frac{3}{2a}-\frac{3}{2b}$ so that $\beta(6, 6)=0$.

As for Lemma 5.1 in \cite{K-T}, their formulation is for $\mathbb C$-valued functions in $L^p$, while we need it for $L^2$ valued functions in $L^p$ (that is, $F\in L^p (dx) L^2(dy)$). We have the following analog:

\begin{lemma} Let $1<p<\infty$. Any $ F\in L^p (dx) L^2(dy)$ can be written as
\begin{align*}
F(x, y)= \sum c_k \chi_k(x, y)
\end{align*}
where each $c_k \ge 0$,  $\|\chi_k(x,y)\|_{L^2(dy)}$ is supported in $x$ in a set of measure $O(2^k)$,  $\|\chi_k \|_{L^{\infty}(dx)L^2(dy)}
\le C 2^{- \frac{k}{p}}$ and $\sum c_k^p \le C \|F\|^p_{L^p(dx) L^2 (dy)}$.
\end{lemma}
\begin{proof}

Define, for $\alpha >0$,
\begin{align*}
\lambda(\alpha)=\big|\big\{ \|F(x, \cdot)\|_{L^2(\mathbb R^3)} > \alpha\big\}\big|
\end{align*}
and
\begin{align*}
&\alpha_k= \inf_{\lambda(\alpha) < 2^k} \alpha\\
&c_k = 2^{\frac{k}{p}} \alpha_k\\
&\mbox{and define}\\
&\chi_k(x, y)=
\begin{cases}
&\frac{1}{c_k} F(x, y) \, \, \mbox{if} \, \, \alpha_{k+1} < \|F(x, \cdot)\|_{L^2(dy)} \le \alpha_k\\
& 0 \, \, \mbox{ otherwise}.
\end{cases}
\end{align*}
From here, we get right away
\begin{align*}
&\|\chi_k(x, \cdot)\|_{L^2}=
\begin{cases}
&\frac{1}{c_k} \|F(x, \cdot)\|_{L^2} \, \, \mbox{if} \, \, \alpha_{k+1} < \|F(x, \cdot)\|_{L^2(dy)} \le \alpha_k\\
& 0 \, \, \mbox{ otherwise}.
\end{cases}
\end{align*}
Thus
\begin{align*}
\|F(x, \cdot)\|_{L^2(dy)} =\sum c_k \|\chi_k(x, \cdot)\|_{L^2(dy)}.
\end{align*}
is exactly the atomic decomposition  of \cite{K-T} corresponding to the $L^p$  function $x \to \|F(x, \cdot)\|_{L^2(dy)}$.
From here we get for free $\|\chi_k(x,y)\|_{L^2(dy)}$ is supported in $x$ in a set of measure $O(2^k)$,  $\|\chi_k \|_{L^{\infty}(dx)L^2(dy)}
\le C 2^{- \frac{k}{p}}$ and $\sum c_k^p \le C \|F\|^p_{L^p(dx) L^2 (dy)}$.
\end{proof}

To finish the proof, following \cite{K-T}, use the above decomposition to write
\begin{align*}
&F(t, x, y)=\sum f_k(t)F_k (t, x, y)\, \, \mbox{ (thus $c_k$ is called $f_k$, $\chi_k$ is called $F_k$)}\\
&G(t, x, y)=\sum g_k(t)G_k (t, x, y)
\end{align*}
thus
\begin{align*}
\sum |T_j(F, G)| \le \sum |T_j( f_k F_k,  g_l G_l)|
\end{align*}
and
optimizing there exists $\epsilon >0$ such that
\begin{align*}
|T_j( f_k F_k,  g_l G_l)| \lesssim 2^{-\epsilon \left(|k - \frac{3}{2} j|+|l - \frac{3}{2} j|\right)}\|f_k\|_{L^2}\|g_l\|_{L^2}
\end{align*}
which can be summed as in \cite{K-T}:
\begin{align*}
&\sum_{j, k, l}|T_j( f_k F_k,  g_l G_l)|
\lesssim \sum_{ k, l}
2^{-\epsilon' \left(|k -l|\right)}\|f_k\|_{L^2}\|g_l\|_{L^2}\\
&\lesssim \left(\sum_{ k} \|f_k\|^2_{L^2}\right)^{\frac{1}{2}}
\left(\sum_{ k} \|g_k\|^2_{L^2}\right)^{\frac{1}{2}}\\
&\lesssim \left(\sum_{ k} \|f_k\|^{\frac{6}{5}}_{L^2}\right)^{\frac{5}{6}}
 \left(\sum_{ k} \|g_k\|^{\frac{6}{5}}_{L^2}\right)^{\frac{5}{6}}\\
 &\lesssim \|F\|_{L^{\frac{6}{5}}(dx) L^2(dy)}\|G\|_{L^{\frac{6}{5}}(dx) L^2(dy)}.
\end{align*}

\section{Proof of Theorem \ref{D^2}}

\subsection{A Priori Bounds and basic estimates}
We will use the following estimates:

 \begin{proposition}\label{basicphi}
 For any smooth, $L^2$, self-adjoint, positive semi-definite kernel $\Gamma(x, y)$  we have the pointwise estimates
 \begin{equation}\label{gammaprop}
    |\Gamma( x,y)|^2\leq \Gamma( x,x)\Gamma( y,y),
\end{equation}
and
 \begin{equation}\label{gammaprop2}
   \big| \nabla_x \Gamma( x,z)\big| \leq E_k( x)^{\frac{1}{2}}\cdot \Gamma( z,z)^{\frac{1}{2}},
\end{equation}
where $E_k(t, x)$ is the kinetic energy density defined as
\begin{align}
E_k( x) = \nabla_{x}\cdot\nabla_{y}\Gamma(t, x, y)\bigg|_{x=y}. \label{kin}
\end{align}
 \end{proposition}
 \begin{proof}
 The above two estimates follow from the Cauchy-Schwarz inequality,  and writing
 \begin{equation}
     \Gamma(x,y)=\sum_{i} \lambda_i \psi_{i}(x)\bar{\psi_i}(y).
 \end{equation}
 \end{proof}
 \begin{proposition}[Fixed time estimates based on conserved quantities]\label{energy1}
 Under the assumptions of Theorem \ref{goal},
 \begin{align*}
 &\|\Gamma(t, x, x)\|_{L^{\infty}(dt)L^1(dx)}= \|\Gamma(0, x, x)\|_{L^1(dx)}=1,\\
 & ||\Gamma(t, x, x)||_{L^{\infty}(dt)L^2(dx )},\\
 & \lesssim||\nabla_x \nabla_y \Gamma||_{L^{\infty}(dt)L^2(dx dy)}
 +\|\Gamma(t, x, x)\|_{L^{\infty}(dt)L^1(dx)}
  \lesssim 1,\\
 & ||\phi||_{L^{\infty}(dt)H^1(dx )} \lesssim 1,\\
 &\|E_k\|_{L^{\infty}(dt)L^1(dx )} \lesssim 1.
\end{align*}
\end{proposition}
 \begin{proposition}[Space-time estimates based on interaction Morawetz]\label{priori}
 Under the assumptions of Theorem \ref{goal},
 \begin{align}
||\Gamma(t,x,x)||_{L^2_{t,x}}  \lesssim 1 \label{morawetzGamma}
\end{align}
which implies
\begin{align}
    ||\phi||_{L^4_tL^4_x} \lesssim 1. \label{morawetz}
\end{align}
\end{proposition}
\begin{proof} A proof of this result has already appeared in the unpublished thesis \cite{Jackythesis}.
For completeness, we include the proof in section \ref{pfMorawetz}. \end{proof}




 \subsection{Estimates for the RHS of  \eqref{leq} in dual Strichartz norms}\label{dual}

 Denote
  \begin{equation}\label{main99}
    \left( \S +\frac{1}{N} V_N(x-y)\right) \Lambda(t,x,y)=\textmd{Term1}+\textmd{Term2}+\textmd{Term3} +\textmd{Term4},
\end{equation}
where
\begin{align*}
&    \textmd{Term1}=-(V_N\ast \Gamma(t,x,x)+V_N\ast \Gamma(t,y,y))\cdot \Lambda(t,x,y),\\
&\textmd{Term2}=V_N\Lambda \circ \Gamma+\bar{\Gamma} \circ V_N\Lambda ,\\
&
    \textmd{Term3}=\Lambda \circ V_N\Gamma+V_N\bar{\Gamma} \circ \Lambda ,\\
&\mbox{and}\\
&\textmd{Term4}=2(V_N\ast|\phi|^2)(y)\phi(x)\phi(y)+2(V_N\ast|\phi|^2)(x)\phi(x)\phi(y).
\end{align*}
Let $2 < p_0 \le \frac{8}{3}$ and define the localized, restricted Strichartz norm
\begin{align*}
&\|\Lambda\|_{\mathcal{S}_{restrited}[T_1, T_2]}\\
&=  \sup_{p_0\le p \le \infty,\, \,  p, q \, \, \text{admissible}}\|\Lambda\|_{L^{p}[T_1, T_2] L^{q}(dx) L^2(dy)}\\
&\quad + \sup_{p_0\le p \le \infty, \, \,  p, q \, \, \text{admissible}}\|\Lambda\|_{L^{p}[T_1, T_2] L^{q}(dy) L^2(dx)} \\
 &\quad +
\sup_{2 \le p \le \infty, \, \,  p, q \, \, \text{admissible}}\|\Lambda\|_{L^p[T_1, T_2] L^q(d(x-y)) L^2(d(x+y))} \}.
\end{align*}
and, for $(p, q)$ an admissible Strichartz pair, define the localized dual norms
\begin{align*}
&\|F\|_{\mathcal{S}_{dual}^{p', q'}[T_1, T_2]}\\
&= \min \big\{ \|F\|_{L^{p'}[T_1, T_2] L^{q'}(dx) L^2(dy)},
 \|F\|_{L^{p'}[T_1, T_2] L^{q'}(dy) L^2(dx)},  \|F\|_{L^{p'}[T_1, T_2] L^{q'}(d(x-y)) L^2(d(x+y))} \}.
\end{align*}

In preparation for applying Theorem \ref{concise},  we state the following estimates, in a simple (but not sharp)  form which will suffice for our goal. We will use Proposition \ref{basicphi}, Proposition \ref{energy1} and Proposition \ref{priori} to bound various terms uniformly in $N$, keeping track only of $\|\Gamma(t, x, x)\|_{L^2([T_1, T_2])}$ which will be small (after suitably localizing in time), and
$\|\Lambda\|_{\mathcal{S}_{restrited}[T_1, T_2]}$ which will be handled by a bootstrapping argument.
\begin{theorem} \label{simple} Under the assumptions of Theorem \ref{goal},
for $k=1, 2, 3$ we have
\begin{align*}
&\|\textmd{Term \, k}\|_{\mathcal{S}_{dual}^{\frac{8}{5}, \frac{4}{3}}[T_1, T_2]}
\lesssim N^{\frac{1}{2}} ||\Gamma(t,x,x)||^{\frac{1}{4}}_{L^2[T_1, T_2]L^2(dx)} \|\Lambda\|_{\mathcal{S}_{restrited}[T_1, T_2]},\\
&\|\nabla \textmd{Term \, k}\|_{\mathcal{S}_{dual}^{\frac{8}{5}, \frac{4}{3}}[T_1, T_2]}
\lesssim  N^{\frac{3}{2}} ||\Gamma(t,x,x)||^{\frac{1}{4}}_{L^2[T_1, T_2]L^2(dx)} \|\Lambda\|_{\mathcal{S}_{restrited}[T_1, T_2]}\\
&+ N^{\frac{1}{2}} ||\Gamma(t,x,x)||^{\frac{1}{4}}_{L^2[T_1, T_2]L^2(dx)} \|\nabla \Lambda\|_{\mathcal{S}_{restrited}[T_1, T_2]},\\
&\|\nabla_x \nabla_y \textmd{Term \, k}\|_{\mathcal{S}_{dual}^{\frac{8}{5}, \frac{4}{3}}[T_1, T_2]}
\lesssim  N^{\frac{3}{2}} ||\Gamma(t,x,x)||^{\frac{1}{4}}_{L^2[T_1, T_2]L^2(dx)} \|\nabla \Lambda\|_{\mathcal{S}_{restrited}[T_1, T_2]}\\
&+||\Gamma(t,x,x)||^{\frac{1}{2}}_{L^2[T_1, T_2]L^2(dx)} \\
    &\cdot \left(||\nabla_x\nabla_y\Lambda(t,x,y)||_{L^{\frac{8}{3}}[T_1, T_2]L^4(dx)L^2(dy)}+
    ||\nabla_x\nabla_y\Lambda(t,x,y)||_{L^{\frac{8}{3}}[T_1, T_2]L^4(dy)L^2(dx)}\right).
\end{align*}

Also,
\begin{align*}
&\|\textmd{Term4}\|_{\mathcal{S}_{dual}^{2, \frac{6}{5}}[T_1, T_2]}
\lesssim  1,\\
&\|\nabla \textmd{Term4}\|_{\mathcal{S}_{dual}^{2, \frac{6}{5}}[T_1, T_2]}
\lesssim   N, \\
&\|\nabla_x \nabla_y \textmd{Term4}\|_{\mathcal{S}_{dual}^{2, \frac{6}{5}}[T_1, T_2]}
\lesssim   N.
\end{align*}

\end{theorem}
Notice that $\nabla_x \nabla_y \textmd{Term4}$  had to be estimated in an end-point dual Strichartz norm.

The proof of this theorem is based on Proposition \ref{basicphi}, Proposition \ref{energy1},
Proposition \ref{priori}
 and H\"older's inequality. It will be given in an appendix.

 \subsection{Polynomial in $N$ estimates for the Strichartz norms of $\Lambda$  and its derivatives.}
In this subsection, we finish the proof of Theorem \ref{D^2}.

Using the a priori estimates of  Theorem \ref{simple}, as well as the Strichartz estimates of Theorem \ref{concise}, we estimate
first $\big\| \Lambda \big\|_{\mathcal{S}_{restricted}}$ and then use this to estimate $\big\| \nabla \Lambda \big\|_{\mathcal{S}_{restricted}}$ and then $\big\|\nabla_x \nabla_y \Lambda \big\|_{\mathcal{S}_{restricted}}$.

 \begin{theorem}\label{maincontrol}
Under the assumptions of Theorem \ref{goal}, the following holds
\begin{align*}
   & \big\| \Lambda \big\|_{\mathcal{S}_{restricted}[0, \infty)}
    \lesssim N^4.
\end{align*}
 \end{theorem}

 \begin{proof}

Recall
 \begin{equation}\label{main9}
    \left( \S +\frac{1}{N} V_N(x-y)\right) \Lambda(t,x,y)=\textmd{Term1}+\textmd{Term2}+\textmd{Term3} +\textmd{Term4}.
\end{equation}

Adapting the argument of Bourgain \cite{Bourgain98}, we use estimate \eqref{morawetzGamma} to break up $[0, \infty)$ into
about $N^4$  time intervals $[T_j, T_{j+1}]$ where
 where
$N^{\frac{1}{2}}||\Gamma(t,x,x)||^{\frac{1}{4}}_{L^2[T_j, T_{j+1}]L^2(dx)} \le \epsilon$ (with  $\epsilon$ sufficiently small to be determined later).

We will show that each $\big\| \Lambda \big\|_{\mathcal S [T_j, T_{j+1}]} \le C$ where $C$ depends only on the initial conditions of the system at $t=0$.

For $t \in [T_j, T_{j+1}]$ we have
\begin{align}
&\Lambda(t) = e^{i t \left(-\Delta_{x, y} +\frac{1}{N} V_N \right)} \Lambda(T_j) + i
\sum_{k=1}^4 \int_{T_j}^t  e^{i (t-s) \left(-\Delta_{x, y} +\frac{1}{N} V_N \right)}  \textmd{Term  k}(s) ds \notag\\
&:= e^{i t \left(-\Delta_{x, y} +\frac{1}{N} V_N \right)} \Lambda(T_j) +
\sum_{k=1}^4 \Lambda_k. \label{sum}
\end{align}
Using Theorem \eqref{concise}, and the conservation \eqref{L^2cons}
\begin{align*}
&|| e^{i t \left(-\Delta_{x, y} +\frac{1}{N} V_N \right)} \Lambda(T_j)||_{\mathcal S [T_j, T_{j+1}]}
\lesssim \|\Lambda(T_j)\|_{L^2} \lesssim 1.
\end{align*}

Also Theorem \ref{concise} and Theorem \ref{simple} imply,
\begin{align*}
&\|\sum_{k=1}^4 \Lambda_k\|_{\mathcal{S}_{restricted}[T_j, T_{j+1}]}\le
\sum_{k=1}^4 \|\Lambda_k\|_{\mathcal{S}_{restricted}[T_j, T_{j+1}]}\\
&\lesssim
 N^{\frac{1}{2}} ||\Gamma(t,x,x)||^{\frac{1}{4}}_{L^2[T_j, T_{j+1}]L^2(dx)} \|\Lambda\|_{\mathcal{S}_{restrited}[T_j, T_{j+1}]} + 1\\
&\lesssim \epsilon \| \Lambda\|_{\mathcal{S}_{restricted}[T_j, T_{j+1}]} +1.
\end{align*}

Putting everything together, using the decomposition \eqref{sum},
\begin{align*}
 &||\Lambda||_{\mathcal{S}_{restricted}[T_j, T_{j+1}]}
 \le C_1 + C_2 \epsilon||\Lambda||_{\mathcal{S}_{restricted}[T_j, T_{j+1}]}
\end{align*}
where $C_1$, $C_2$ depend only on the initial conditions of the system at time $t=0$. If we choose  $C_2 \epsilon < \frac{1}{2}$, we get
\begin{align}
 &||\Lambda||_{\mathcal{S}_{restricted}[T_j, T_{j+1}]}\le 2 C_1 \label{noder} \\
 &\mbox{ and, summing over all $\sim N^4$ intervals,}\notag\\
 &||\Lambda||_{\mathcal S[0, \infty)} \lesssim N^4. \notag
\end{align}
\end{proof}

 \begin{theorem}\label{onederiv}
Under the assumptions of Theorem \ref{goal}, the following holds
\begin{align*}
   & \big\|\nabla  \Lambda \big\|_{{\mathcal S}_{restricted}}
    \lesssim N^{5}.
\end{align*}
 \end{theorem}
 \begin{proof}
The proof uses the estimates of Theorem \eqref{maincontrol}, and is similar in structure. It uses the same $\sim N^4$ intervals $[T_j, T_{j+1}]$.

Differentiate the equation \eqref{main9}, and estimate the right-hand side in a dual Strichartz space.

Thus
 \begin{align*}
   & \left( \S +\frac{1}{N} V_N(x-y)\right) \nabla \Lambda(t,x,y)\\
   &= \nabla\textmd{Term1}+ \nabla\textmd{Term2}+ \nabla\textmd{Term3} + \nabla\textmd{Term4}\\
    &-  \nabla\left(\frac{1}{N} V_N(x-y)\right)\Lambda.
\end{align*}
Call the last term $\textmd{Term5}$.
Following the argument of the previous proof:

\begin{align*}
&\nabla \Lambda(t) = e^{i t \left(-\Delta_{x, y} +\frac{1}{N} V_N \right)} \nabla \Lambda(T_j)\\
& + i
\sum_{k=1}^4 \int_{T_j}^t  e^{i (t-s) \left(-\Delta_{x, y} +\frac{1}{N} V_N \right)} \nabla \textmd{Term  k}(s) ds
+\int_{T_j}^t  e^{i (t-s) \left(-\Delta_{x, y} +\frac{1}{N} V_N \right)}  \textmd{Term  5}(s) ds
\notag\\
&:= e^{i t \left(-\Delta_{x, y} +\frac{1}{N} V_N \right)} \nabla \Lambda(T_j) +
\sum_{k=1}^5 \Lambda_k.
\end{align*}

Using conservation of energy (see \eqref{lambdaest}), we have
\begin{align*}
\| e^{i t \left(-\Delta_{x, y} +\frac{1}{N} V_N \right)} \nabla \Lambda(T_j)\|_{\mathcal{S}_{restricted}[T_j, T_{j+1}]} \lesssim
\| \nabla \Lambda(T_j)\|_{L^2} \lesssim 1.
\end{align*}

It remains to estimate $ \nabla \textmd{Term  1}, \cdots,  \nabla \textmd{Term  4}$ and $  \textmd{Term  5}$ in $\mathcal{S}_{dual}^{p', q'}$.

We have,
using H\"older's inequality
\begin{align*}
&\|\textmd{Term5}\|_{L^2[T_j, T_{j+1}]L^{\frac{6}{5}}(d(x-y))L^2(d(x+y))} \lesssim N \|\Lambda\|_{L^2[T_j, T_{j+1}] L^6(d(x-y))L^2(d(x+y))}\\
&\lesssim N\, \, \mbox{ (we used  \eqref{noder})},
\end{align*}
while, from Theorem \ref{concise} and Theorem \ref{simple} and another application of \eqref{noder},
\begin{align*}
&\sum_{k=1}^4 \|\Lambda_k\|_{\mathcal{S}_{restricted}[T_j, T_{j+1}]}\\
& \lesssim \sum_{k=1}^3 \|\nabla  \textmd{Term \, k}\|_{\mathcal{S}_{dual}^{\frac{8}{5}, \frac{4}{3}}[T_1, T_2]}
+\|\nabla  \textmd{Term 4}\|_{\mathcal{S}_{dual}^{2, \frac{6}{5}}[T_1, T_2]}
\\
&\lesssim N^{\frac{1}{2}}||\Gamma(t,x,x)||^{\frac{1}{4}}_{L^2[T_1, T_2]L^2(dx)} \| \nabla \Lambda\|_{\mathcal{S}_{restrited}[T_1, T_2]}\\
&+ N^{\frac{3}{2}} ||\Gamma(t,x,x)||^{\frac{1}{4}}_{L^2[T_1, T_2]L^2(dx)} \|\Lambda\|_{\mathcal{S}_{restrited}[T_1, T_2]}
+ N\\
&\le C_1 N + C_2 \epsilon || \nabla \Lambda\|_{\mathcal{S}_{restrited}[T_1, T_2]}.
\end{align*}
Since $\epsilon$ is chosen so that $C_2 \epsilon < \frac{1}{2}$, summing the previous estimates we get $\|\nabla \Lambda\|_{\mathcal{S}_{restricted}[T_j, T_{j+1}]} \lesssim N$ and, summing over all $\sim N^4$ intervals,
\begin{align*}
   & \big\|\nabla  \Lambda \big\|_{{\mathcal S}_{restricted}}
    \lesssim N^{5}.
\end{align*}

\end{proof}

Finally,
\begin{theorem}\label{twoderiv}
Under the assumptions of Theorem \ref{goal}, the following holds
\begin{align*}
   & \big\|\nabla_x \nabla_y  \Lambda \big\|_{{\mathcal S}_{restricted}}
    \lesssim N^{\frac{13}{2}}.
\end{align*}
 \end{theorem}
 \begin{proof}
We write
\begin{align*}
   & \left( \S +\frac{1}{N} V_N(x-y)\right) \nabla_x\nabla_y \Lambda(t,x,y)\\
   &= \nabla_x\nabla_y\textmd{Term1}+\cdots + \nabla_x\nabla_y\textmd{Term4}\\
    &-  \nabla_x\left(\frac{1}{N} V_N(x-y)\right)\nabla_y\Lambda-  \nabla_y\left(\frac{1}{N} V_N(x-y)\right)\nabla_x\Lambda
    -\nabla_x\nabla_x\left(\frac{1}{N} V_N(x-y)\right)\Lambda
\end{align*}
with initial conditions $\|\nabla_x\nabla_y \Lambda_0\|_{ L^2} \lesssim N$.
Unlike the previous two proofs, we no longer have a priori bounds on the growth of $\|\nabla_x\nabla_y \Lambda(t)\|_{ L^2}$ - in fact this is what we are trying to prove. Now we split $[0, \infty)$ differently than before. Now we only require $||\Gamma(t,x,x)||^{\frac{1}{2}}_{L^2[T_j, T_{j+1}]L^2(dx)} \le \epsilon$, with $\epsilon $ (independent of $N$) to be determined later. The number of intervals only depends on $||\Gamma(t,x,x)||_{L^2[0, \infty)L^2(dx)} \lesssim 1$, and is independent of $N$.
 We apply
  Theorem \ref{concise} and Theorem \ref{simple} directly on $[T_i, T_{i+1}]$, using the estimates for $\big\| \Lambda \big\|_{{\mathcal S}_{restricted}}$ and $\big\| \nabla\Lambda \big\|_{{\mathcal S}_{restricted}}$ from the previous two theorems.

  For $k=1, 2, 3$ we have
  \begin{align*}
  &\|\nabla_x \nabla_y \textmd{Term \, k}\|_{\mathcal{S}_{dual}^{\frac{8}{5}, \frac{4}{3}}[T_i, T_{i+1}]}
\lesssim  N^{\frac{3}{2}} ||\Gamma(t,x,x)||^{\frac{1}{4}}_{L^2[T_i, T_{i+1}]L^2(dx)} \|\nabla \Lambda\|_{\mathcal{S}_{restrited}[T_i, T_{i+1}]}\\
&+||\Gamma(t,x,x)||^{\frac{1}{2}}_{L^2[T_i, T_{i+1}]L^2(dx)} \|\nabla_x\nabla_y \Lambda\|_{\mathcal{S}_{restrited}[T_i, T_{i+1}]} \\
&\le C_1
  N^{\frac{3}{2}} N^5 + C_2 \epsilon \|\nabla_x\nabla_y \Lambda\|_{\mathcal{S}_{restrited}[T_i, T_{i+1}]},
\end{align*}
while
\begin{align*}
  &\|\nabla_x \nabla_y \textmd{Term \, 4}\|_{\mathcal{S}_{dual}^{2, \frac{6}{5}}[T_i, T_{i+1}]}\lesssim N.
  \end{align*}
  As for the terms where the derivatives fall on the potential, for example
  \begin{align*}
 &\| \nabla_x\nabla_x\left(\frac{1}{N} V_N(x-y)\right)\Lambda\| _{\mathcal{S}_{dual}^{2, \frac{6}{5}}[T_i, T_{i+1}]}
 \lesssim \|\nabla^2 \frac{1}{N} V_N\|_{L^{\frac{3}{2}}}\|\Lambda\| _{L^2[T_i, T_{i+1}]L^6(d(x-y))L^2(d(x+y))}\\
 &\lesssim N^2 \|\Lambda\|_{\mathcal{S}_{restrited}[T_i, T_{i+1}]} \lesssim N^6.
  \end{align*}
  Thus, with some choice of constants $C_i$ depending only on the initail conditions,  we get from Theorem \ref{concise}
  \begin{align*}
  & \|\nabla_x\nabla_y \Lambda\|_{\mathcal{S}_{restrited}[T_i, T_{i+1}]}\\
   & \le C_1\|\nabla_x\nabla_y \Lambda(T_i)\|_{L^2} + C_2 N^{\frac{13}{2}} + C_3 \epsilon \|\nabla_x\nabla_y \Lambda\|_{\mathcal{S}_{restrited}[T_i, T_{i+1}]}.
   \end{align*}
   If we pick $ C_3 \epsilon < \frac{1}{2}$,
   and notice $\|\nabla_x\nabla_y \Lambda(T_i)\|_{L^2} \le \|\nabla_x\nabla_y \Lambda\|_{\mathcal{S}_{restrited}[T_{i-1}, T_i]}$,
    we conclude
    \begin{align*}
  & \|\nabla_x\nabla_y \Lambda\|_{\mathcal{S}_{restrited}[T_i, T_{i+1}]}\\
   & \le 2 \left( C_1\|\nabla_x\nabla_y \Lambda\|_{\mathcal{S}_{restrited}[T_{i-1}, T_i]} + C_2 N^{\frac{13}{2}} \right).
   \end{align*}

  Applying this arguments a finite number of times (independent of $N$), and summing the result, we are done.

    \end{proof}

    \section{Proof of Proposition \ref{priori} \label{pfMorawetz}}

    The outline of this section is inspired, in part, by \cite{CKSTT},  \cite{CPTz} and the similarities between the HFB system and the GP hierarchy. The main result appeared in
    the unpublished thesis \cite{Jackythesis}.

\subsection{Local Conservation Laws } Let us start by defining the relevant quantities which will allow us to effectively capture the conservation laws of the HFB system. We define
\begin{subequations}
\begin{align}
&T_{00}=\rho :=\ \Gamma(x; x) \\
&T_{j0}=T_{0j}=P_j:=\ \frac{1}{2i} \int dx'\ \delta(x-x') [\bd_{x_j'} \Gamma(x; x')-\bd_{x_j}\Gamma(x; x')]\\
&T_{jk}=\sigma_{jk}+p\delta_{jk}:=\ \int dx'\ \delta(x-x')(\bd_{x_j}\bd_{x_k'}+\bd_{x_k}\bd_{x_j'})\Gamma(x; x')\\
&\qquad\qquad\qquad\qquad\ +\delta_{jk}\frac{1}{2}\left(-\lapl \rho+\int dy\ V_N(x-y)\mathcal{L}(x, y; x, y) \right) \nonumber\\
&l_j = \frac{1}{2}\int dy\ V_N(x-y)\{\bd_{y_j}\mathcal{L}(x, y; x, y)-\bd_{x_j}\mathcal{L}(x, y; x, y)\}\\
&\mathcal{L}(x, y; x', y') :=\  \Gamma(x; x')\Gamma(y; y')+\Gamma(x; y')\Gamma(y; x')\\
&\qquad\qquad\qquad\quad  +\overline{\Lambda}(x, y)\Lambda(x', y')-2\bar \phi(x)\bar\phi(y)\phi(x')\phi(y'). \nonumber
\end{align}
\end{subequations}
In the literature, $T_{\mu\nu}$ is often referred to as the pseudo-stress-energy tensor and $L$ is the two-particle marginal density matrix of our quasifree state.  Then the associated
local conservation laws are given by
\begin{equation}\label{local-conser-law}
\begin{cases}
\bd_t\rho +2\grad \cdot P = 0\\
\bd_tP+\grad\cdot (\boldsymbol{\sigma} +p \vect{I})+ l=0
\end{cases}.
\end{equation}

To derive the local conservation laws, it is convenient to first rewrite the equation for $\Gamma(x; x')$ in the following form
\begin{align}\label{bbgky-gamma}
\left\{\frac{1}{i}\frac{\bd}{\bd t} + \lapl_{x}-\lapl_{x'}\right\} \Gamma(x; x') = B_{V}(\mathcal{L})
\end{align}
where
\begin{subequations}
\begin{align}
B_{V}(\mathcal{L}) :=&\ B_{V}^+(\mathcal{L})-B_{V}^-(\mathcal{L}),\\
B_{V}^+(\mathcal{L})(x; x'):=&\ \int dydy'\ V_N(x-y)\delta(y-y') \mathcal{L}(x, y; x', y'),\\
B_{V}^-(\mathcal{L})(x; x'):=&\ \int dydy'\ V_N(x'-y) \delta(y-y')\mathcal{L}(x, y; x', y').
\end{align}
\end{subequations}
Notice \eqref{bbgky-gamma} has the structure of a BBGKY hierarchy, that is, the evolution of the lower marginal density matrix depends on the higher marginal density. Unlike, the standard BBGKY
hierarchy, the quasifree structure of our state allows us to decompose our two-particle marginal density matrix $\mathcal{L}$ into a linear combination of products of one-particle marginal densities $\Gamma, \Lambda$ and the condensate wave function $\phi$.

\begin{proposition}
Let $\Gamma$ be a  smooth solution to \eqref{bbgky-gamma}, then we have the local conservation of number
\begin{align}\label{conserved-num}
\frac{\bd \rho}{\bd t} + 2\grad\cdot P = 0.
\end{align}
\end{proposition}

\begin{proof}
By direct calculation, we see that
\begin{subequations}
\begin{align}
\bd_t\rho =&\  \int \frac{du du'}{(2\pi)^6}\ e^{i(u-u')\cdot x}\bd_t \widehat \Gamma(u; u') \nonumber\\
=&\ i \int \frac{du du'}{(2\pi)^6}\ e^{i(u-u')\cdot x}(u^2- (u')^2)\widehat \Gamma(u; u') \label{number-term1}\\
&\ +i \int  \frac{du du'}{(2\pi)^6}\ e^{i(u-u')\cdot x}\widehat {B_{V}(\mathcal{L})}(u; u'). \label{number-term2}
\end{align}
\end{subequations}
For the first term, we have that
\begin{align*}
\eqref{number-term1} =  \grad_x\cdot \int \frac{du du'}{(2\pi)^6}\ e^{i(u-u')\cdot x}(u+ u')\widehat \Gamma(u; u') = -2\grad_x\cdot P.
\end{align*}
For the second term,  we have that  $\eqref{number-term2}= iB_{V}(\mathcal{L})(x; x)=0$.
\end{proof}

\begin{proposition}
Let $(\phi, \Gamma, \Lambda)$ be a  smooth solution to the HFB system, then we have the continuity equation
\begin{align}\label{conserved-momentum}
 \bd_tP+\grad\cdot (\boldsymbol{\sigma} +p \vect{I})+ l=0.
\end{align}
\end{proposition}

\begin{proof}
Differentiating $P$ with respect to time yields
\begin{align*}
\bd_tP(x) =&\ \frac{1}{i}\int \frac{du du'}{(2\pi)^6}\ e^{i(u-u')\cdot x}\frac{(u+u')}{2}(u^2-(u')^2)\widehat \Gamma(u; u')\\
&\ +\frac{1}{i}\int \frac{du du'}{(2\pi)^6}\ e^{i(u-u')\cdot x} \frac{(u+u')}{2}\widehat{B_{V}(\mathcal{L})}(u; u')\\
=&\ -\frac{1}{2}\grad_x\cdot\int \frac{du du'}{(2\pi)^6}\ e^{i(u-u')\cdot x}(u+u')\otimes (u+u')\widehat \Gamma(u; u')\\
&\ +\frac{1}{i}\int \frac{du du'}{(2\pi)^6}\ e^{i(u-u')\cdot x} \frac{(u+u')}{2}\widehat{B_{V}(\mathcal{L})}(u; u')=:\ J_1+J_2.
\end{align*}
Let us first handle the $J_1$ term. Notice we have that
\begin{align*}
J_1  =&\  -\frac{1}{2}\grad_x\cdot \int \frac{dudu'}{(2\pi)^6}\ e^{i(u-u')\cdot x} (u-u')^{\otimes 2}\widehat{\Gamma}(u; u')\\
&\ - \grad_x\cdot \int \frac{dudu'}{(2\pi)^6}\ e^{i(u-u')\cdot x} (u\otimes u' + u'\otimes u) \widehat{\Gamma}(u; u').
\end{align*}
Then, completing the Fourier inversion gives us
\begin{align*}
J_1=&\ \frac{1}{2}\grad\cdot \grad^2\rho(x)-\grad\cdot \boldsymbol{\sigma} = -\grad\cdot \left(-\frac{1}{2}\lapl \rho \vect{I}+\boldsymbol{\sigma}\right).
\end{align*}

Next, we deal with the $J_2$ term. By the Fourier inversion, we write
\begin{align*}
J_2 =&\ -\frac{1}{2}\int dx'\ \delta(x-x') \left\{ \grad_x B_V(\mathcal{L})(x; x')- \grad_{x'} B_V(\mathcal{L})(x; x')\right\}.
\end{align*}
Then we observe that
\begin{align*}
&\int dx'\ \delta(x-x') \grad_x B_V(\mathcal{L})(x; x')\\
&=\int dx' dy\ \delta(x-x')\grad_x\left(\{V_N(x-y)-V_N(x'-y)\}\mathcal{L}(x, y; x', y)\right)\\
&=\int dx' dy\ \delta(x-x')\grad_x\left(V_N(x-z)\right)\mathcal{L}(x, y; x', y)\\
&\quad\ +\int dx' dy\ \delta(x-x')\{V_N(x-y)-V_N(x'-y)\}\grad_x\mathcal{L}(x, y; x', y)\\
& = \int dy\ \grad_x\left(V_N(x-y)\right)\mathcal{L}(x, y; x, y).
\end{align*}
Likewise, we have that
\begin{align*}
&\int dx'\ \delta(x-x') \grad_{x'} B_V(\mathcal{L})(x; x') = -\int dy\ \grad_x\left(V_N(x-y)\right)\mathcal{L}(x, y; x, y).
\end{align*}
Hence it follows
\begin{align*}
J_2 &= -\int dy\ \grad_x\left(V_N(x-y)\right)\mathcal{L}(x, y; x, y)\\
& = \frac{1}{2}\int dy\ \{\grad_yV_N(x-y)-\grad_xV_N(x-y)\}\mathcal{L}(x, y; x, y)\\
&=-\frac{1}{2}\grad_x\left( \int dy\ V_N(x-y)\mathcal{L}(x, y; x, y)\right) - l\\
&= -\frac{1}{2}\grad_x\cdot \left( \int dy\ V_N(x-y)\mathcal{L}(x, y; x, y)\vect{I}\right) - l.
\end{align*}
This completes the argument.
\end{proof}

\subsection{Interaction Morawetz Estimate} The main result of this section is the interaction Morawetz-type estimate for the $\Gamma$ equation.
To prove the estimate, we need a two-particle Morawetz identity  for the truncated two-particle marginal density matrix
\begin{align}
L(x, y; x', y') = \Gamma(x; x')\Gamma(y; y').
\end{align}

We formally\footnote{In general, we are not certain whether \eqref{virial-potential} and \eqref{morawetz-action} are well-defined. However, since we are interested when $a(x) = |x|$, it can be shown that \eqref{morawetz-action} is well-defined. More precisely, since $\grad a$ is uniformly bounded, then it follows $|M_a(t)| \le C\llp{\rho}_{L^1(dx)}\llp{P}_{L^1(dx)}$ is uniformly bounded for all time.} define the virial interaction potential for $L$ associated to $a \in C(\rr^3)$ by
\begin{align}\label{virial-potential}
V^a(t):= \int dxdy\ a(x-y)L(t, x, y; x, y)
\end{align}
and its corresponding Morawetz action
\begin{align}\label{morawetz-action}
M^a(t):= \bd_tV^a(t)= 2\int dxdy\ \grad a(x-y)\cdot [P(x)\rho(y)-\rho(x)P(y)].
\end{align}
Then we have the following truncated two-particle Morawetz identity.
\begin{proposition}\label{second-Morawetz-id}
Let $(\phi, \Gamma, \Lambda)$ be a smooth solution to the HFB system
 with ${\rm tr} \Gamma(t)=1$ and $E(t) \le C$ (see \eqref{b5-nbr0}, \eqref{energy}),
 and let $a(x) = |x|$. Then we have
the identity
\begin{subequations}\label{second-mora-id}
\begin{align}
\dot M^a(t) =&\ 2 \int dxdy\ (-\lapl \lapl a) (x-y)\rho(x)\rho(y) \label{morawetz2-term1}\\
&\ +\int dxdy\ \lapl a(x-y)\Big\{\rho(x) \int dz\ V_N(y-z)\mathcal{L}(y, z; y, z) \nonumber \\
&\  +\rho(y) \int dz\ V_N(x-z)\mathcal{L}(x, z; x, z) \Big\}  \label{morawetz2-term2}\\
&\ + 2\int dxdy\ \grad^2 a(x-y):\Big\{\boldsymbol{\sigma}(x)\rho(y)+\rho(x)\boldsymbol{\sigma}(y) \nonumber\\
&\ -4P(x)\otimes P(y)\Big\}  \label{morawetz2-term3}\\
&\ +2\int dxdy\ \grad a(x-y)\cdot \left\{\rho(x)l(y)-l(x)\rho(y)\right\}.  \label{morawetz2-term4}
\end{align}
\end{subequations}
Here, $:$ denotes the standard double dot product, that is, for any $n\times n$ matrices $A$ and $B$, we have that $A:B= \sum_{i, j} a_{ij}b_{ij}$.
\end{proposition}

\begin{remark}
Let us note that Proposition \ref{second-Morawetz-id} only states that for each fixed $N$, identity \eqref{second-mora-id} holds. It does not say that the identity is independent of $N$. In fact, we are not sure whether
\eqref{morawetz2-term4} stays uniformly bounded in $N$. However, this does not pose any issues for us since shortly we will see that the term gives a positive contribution which we can ignore when proving the interaction Morawetz estimate.
\end{remark}

\begin{proof}
The main issue is to show that any integration by parts is justified by the conservation laws.  It is convenient to first note some facts about the pseudo stress-energy tensor. By the conservation laws, we see that
$\rho(x) \in L^1(dx)\cap L^3(dx)$, the components of $P(x)$  are in $L^1(dx)\cap L^\frac{3}{2}(dx)$ and the components of $\boldsymbol{\sigma}(x)$ are in $L^1(dx)$.  However, we don't know anything about the decay properties of $\Delta \rho$ appearing in $T_{jk}$.

To handle any issues with the integration by parts, we apply a smooth spatial cutoff function. Let $\chi \in C^\infty_0(\rr^d)$ be a radial function whose support is contained in the ball $B(0, 2)$ and is identically $1$ on $B(0, 1)$. For every $L>0$, define
\begin{align}\label{cutoff-Morawetz-action}
M^a_L(t) := 2\int dxdy\ \chi(\frac{|x-y|}{L})\grad a(x-y)\cdot [P(x)\rho(y)-\rho(x)P(y)].
\end{align}
Taking the time derivative of \eqref{cutoff-Morawetz-action}, applying the local conservation laws \eqref{local-conser-law}, and integrating by parts yields
\begin{subequations}
\begin{align}
\dot M^a_L(t)  =&\  2\int dxdy\ \grad_x\left(\chi(\frac{|x-y|}{L})\grad a(x-y)\right)\nonumber\\
&\ : \Bigg\{\left(-\frac{1}{2}\lapl_x\rho(x)\rho(y)-\rho(x)\frac{1}{2}\lapl_y\rho(y)\right)\vect{I}\label{trunc-morawetz2-term1}\\
&\ + \Bigg(\frac{1}{2}\int dz\ V_N(y-z)\rho(x) \mathcal{L}(y, z; y, z) \label{trunc-morawetz2-term2} \\
&\ + \frac{1}{2}\int dz\ V_N(x-z)\rho(y)\mathcal{L}(x, z; x, z) \Bigg)\vect{I}\nonumber\\
&\ +   \Big\{\boldsymbol{\sigma}(x)\rho(y)+\rho(x)\boldsymbol{\sigma}(y)-4P(x)\otimes P(y)\Big\}\Bigg\} \label{trunc-morawetz2-term3}\\
&\ + 2\int dxdy\ \chi(\frac{|x-y|}{L})\grad a(x-y)\cdot \left\{\rho(x)l(y)-l(x)\rho(y)\right\}.\label{trunc-morawetz2-term4}
\end{align}
\end{subequations}
Next, we consider the limit as $L$ tends to infinity. It is not hard to see that any derivative of $\chi$ is uniformly bounded in $L$ and vanishes near the origin.
Let us first handle \eqref{trunc-morawetz2-term2}. By direct calculation, we have that
\begin{align*}
 &\grad_x\left(\chi(\frac{|x-y|}{L})\grad a(x-y)\right)\\
 & = \frac{1}{L}\chi'(\frac{|x-y|}{L})\frac{(x-y)\otimes(x-y)}{|x-y|^2}+\chi(\frac{|x-y|}{L})\grad^2 a(x-y)
\end{align*}
which means
\begin{subequations}
\begin{align}
 &\eqref{trunc-morawetz2-term2}\nonumber\\
 &=\ \frac{1}{L}\int dxdydz\ \chi'(\frac{|x-y|}{L}) V_N(y-z)\rho(x) \mathcal{L}(y, z; y, z) \label{trunc-error-term1}\\
 &\quad\ + \int dxdydz\ \chi(\frac{|x-y|}{L}) \lapl a(x-y)V_N(y-z)\rho(x) \mathcal{L}(y, z; y, z) \label{trunc-error-term2}\\
 &\quad\ +\text{ similar terms with $x$ and $y$ switched}. \label{trunc-error-term3}
\end{align}

\end{subequations}
Note that by the  conservation of number and energy, we have that
\begin{align*}
|\eqref{trunc-error-term1}| \le \frac{\llp{\chi'}_\infty}{L}\llp{\rho}_{L^1(dx)}\left(\int dydz\  V_N(y-z) \mathcal{L}(y, z; y, z) \right)\rightarrow 0
\end{align*}
as $L\rightarrow \infty$. Next, by the dominated convergence theorem, we see that
$\eqref{trunc-error-term2}+\eqref{trunc-error-term3}\rightarrow \eqref{morawetz2-term2}$.

The term \eqref{trunc-morawetz2-term3} is handled in a similar manner. More precisely, we see that
\begin{subequations}
\begin{align}
\eqref{trunc-morawetz2-term3}=&\ \frac{2}{L}\int dxdy\ \chi'(\frac{|x-y|}{L})\frac{(x-y)\otimes(x-y)}{|x-y|^2} \label{trunc-error-term5}\\
&\ : \Big\{\boldsymbol{\sigma}(x)\rho(y)+\rho(x)\boldsymbol{\sigma}(y)-4P(x)\otimes P(y)\Big\} \nonumber\\
&\ +2\int dxdy\ \chi(\frac{|x-y|}{L})\grad^2 a(x-y)\label{trunc-error-term6}\\
&\ : \Big\{\boldsymbol{\sigma}(x)\rho(y)+\rho(x)\boldsymbol{\sigma}(y)-4P(x)\otimes P(y)\Big\}. \nonumber
\end{align}
\end{subequations}
For the term \eqref{trunc-error-term5}, we have the estimate
\begin{align*}
|\eqref{trunc-error-term5}| \leq&\ \frac{C\llp{\chi'}_\infty}{L} \left(\llp{\rho}_{L^1(dx)} \llp{\boldsymbol{\sigma}}_{L^1(dx)}+ \llp{P}_{L^1(dx)}^2\right)\rightarrow 0
\end{align*}
as $L$ tends to infinity.

For the term  \eqref{trunc-error-term6}, we first recall that $\grad^2 a(x) = |x|^{-1}\left( \vect{I} -\frac{x\otimes x}{|x|^2}\right)$. Then, by Hardy-Littlewood-Sobolev inequality, it follows that
\begin{align*}
|\eqref{trunc-error-term6}| \le&\ \sum_{i, j} \int dx\ \left\{|\sigma_{ij}(x)| (|\cdot|^{-1}\ast \rho)(x)+\int dy\ \frac{|P_i(x)||P_j(y)|}{|x-y|}\right\}\\
 \le&\ C\llp{\boldsymbol{\sigma}}_{L^1(dx)}\llp{|\cdot|^{-1}\ast \rho}_{L^\infty(dx)}+C\llp{P}_{L^\frac{6}{5}(dx)}^2.
\end{align*}
Hence it suffices to check that $(|\cdot|^{-1}\ast \rho)(x)$ is uniformly bounded. Note that we have the estimate
\begin{align*}
\left|\int dy\ \frac{\rho(y)}{|x-y|}\right| \le&\  \left|\int_{|x-y|<1} dy\ \frac{\rho(y)}{|x-y|}\right| +\int dy\ \rho(y)\\
\le&\ \llp{|\cdot|^{-1}}_{L^\frac{3}{2}(B_1(0))} \llp{\rho}_{L^3(dy)}+\llp{\rho}_{L^1(dy)} \le C
\end{align*}
which holds uniformly in $x$. Then, by dominated convergence theorem, we again see that $\eqref{trunc-error-term6}\rightarrow \eqref{morawetz2-term3}$.

Next, for each fixed $N$, we show that $\eqref{trunc-morawetz2-term4} \rightarrow \eqref{morawetz2-term4}$ follows immediately from the Lebesgue dominated convergence theorem. More precisely, we see that
\begin{align*}
|\eqref{trunc-morawetz2-term4}| \le&\ C\int dx\ \rho(x)\int dy\ |l(y)|\\
\le&\ C\llp{\rho}_{L^1(dx)}\left( \llp{V}_{L^1(dx)}\llp{\grad \rho}_{L^{3/2}(dx)}\llp{\rho}_{L^3(dx)}+\llp{V_N}_{L^{3}(dx)}\llp{\grad \Lambda}_{L^2(dxdy)}^2\right)\\
\le&\ CN^{2\beta}.
\end{align*}

Lastly, let us handle \eqref{trunc-morawetz2-term1}. It suffices to estimate
\begin{subequations}
\begin{align}
&\int dxdy\ \grad_x\left(\chi(\frac{|x-y|}{L})\grad a(x-y)\right) : \lapl_x\rho(x)\rho(y)\vect{I} \nonumber\\
&=  \frac{1}{L}\int dxdy\ \lapl_x\chi'(\frac{|x-y|}{L})\rho(x)\rho(y)\label{error-term1}\\
&\quad\ +\int dxdy\ \lapl_x\chi(\frac{|x-y|}{L})\lapl a(x-y) \rho(x)\rho(y)  \label{error-term2}\\
&\quad\ +2\int dxdy\ \grad_x\chi(\frac{|x-y|}{L})\grad_x\lapl a(x-y) \rho(x)\rho(y) \label{error-term3}\\
&\quad\ +\int dxdy\ \chi(\frac{|x-y|}{L}) \lapl\lapl a(x-y) \rho(x)\rho(y).  \label{error-term4}
\end{align}
\end{subequations}
By the remark in the beginning of the proof, we see that
\begin{align*}
|\eqref{error-term1}+\eqref{error-term2}+\eqref{error-term3}| \leq&\ \frac{C}{L}\llp{\rho}_{L^1(dx)}^2
\end{align*}
which converges to zero as $L$ tends to infinity.
Lastly, we have that
\begin{align*}
|\eqref{error-term4}|= 8\pi \int dxdy\ \chi(\frac{|x-y|}{L}) \delta(x-y) \rho(x)\rho(y) = 8\pi \llp{\rho}_{L^1(dx)}^2
\end{align*}
which is clearly uniformly bounded in $L$. Hence, by the dominated convergence theorem, we have the desired result.
\end{proof}

With this special choice of observable, we have that $(-\lapl \lapl a)(x) = 8\pi \delta(x)$ which we have already used.
Also, it is not hard to see that \eqref{morawetz2-term1} and \eqref{morawetz2-term2} are positive terms since
\begin{align}
\int dz\ V_N(x-z)\mathcal{L}(x, z; x, z) \ge 0
\end{align}
given $V_N\ge 0$.  To prove the Morawetz estimate, we need to be able to control \eqref{morawetz2-term3} and \eqref{morawetz2-term4}. In fact, we will show
that $\eqref{morawetz2-term3}\ge 0$ and $\eqref{morawetz2-term2} + \eqref{morawetz2-term4}\ge 0$, then deduce
\begin{align}
 8\pi \int dx\ \rho(t, x)^2\le \bd_tM^a(t)
\end{align}
which will lead to the desired estimate.

\begin{lemma}\label{mora-lem1}
Assume $V_N$ is a positive radial function, i.e. $V_N(x) = N^{3\beta}V(N^\beta |x|)\ge 0$, with $V'(r)\le 0$. Let $(\phi, \Gamma, \Lambda)$
be a  smooth solution to the HFB system.
 Then we have that
$\eqref{morawetz2-term2}+\eqref{morawetz2-term4} \ge 0$.
\end{lemma}

\begin{proof}
By change of variables and integration by parts, we see that
\begin{subequations}
\begin{align}
\eqref{morawetz2-term4} =&\ -4 \int dxdy\ \rho(y)\frac{x-y}{|x-y|} \cdot l(x)\nonumber\\
=&\ -4 \int dxdydz\ N^{4\beta}V'(N^\beta|x-z|)\rho(y) \frac{x-z}{|x-z|}\cdot\frac{x-y}{|x-y|}\mathcal{L}(x, z; x, z) \label{positive-error-term1}\\
&\ - 4 \int dxdydz\ V_N(x-z)\frac{\rho(y)\mathcal{L}(x, z; x, z)}{|x-y|}. \label{positive-error-term2}
\end{align}
\end{subequations}
Notice that $\eqref{positive-error-term2} = -\eqref{morawetz2-term2}$. Finally, exploiting the symmetry $\mathcal{L}(x, z; x, z) = \mathcal{L}(z, x; z, x)$, we can rewrite \eqref{positive-error-term1} as follows
\begin{align*}
\eqref{positive-error-term1} =&\ -2 \int dxdydz\ N^{4\beta}V'(N^\beta|x-z|)\rho(y) \\
&\ \times \left\{\frac{x-z}{|x-z|}\cdot\frac{x-y}{|x-y|}+\frac{z-x}{|z-x|}\cdot\frac{z-y}{|z-y|}\right\}\mathcal{L}(x, z; x, z)\ge 0.
\end{align*}
The last inequality follows from $\mathcal{L}(x, y; x, y)\ge 0$,  $V'(r)\le 0$,  and the identity
\begin{align}
\frac{u-v}{|u-v|}\cdot\frac{u}{|u|}+\frac{v-u}{|v-u|}\cdot\frac{v}{|v|}  = \frac{(|u|+|v|)(1-\cos\theta)}{|u-v|} \ge 0.
\end{align}
\end{proof}

\begin{lemma}\label{mora-lem2}
Let $(\phi, \Gamma, \Lambda)$
be a  smooth solution to the HFB system. Then we have that
$\eqref{morawetz2-term3} \ge 0$.
\end{lemma}

\begin{proof}
Since $\vect{A}(x, y):=\grad^2 a(x-y)$ is symmetric (in fact, it is positive semi-definite), we can rewrite \eqref{morawetz2-term3} by swapping some indices as follows
\begin{subequations}
\begin{align*}
\frac{1}{2}\eqref{morawetz2-term3}=&\ \int dxdydx'dy'\ \delta(x-x')\delta(y-y') \sum_{jk} \bd_{jk}a(x-y) \nonumber\\
&\ \times \Big\{ (\bd_{x_j} \bd_{x'_k}+\bd_{x_k} \bd_{x'_j})+(\bd_{y_j} \bd_{y'_k}+\bd_{y_k} \bd_{y'_j})\nonumber\\
&\quad\ + (\bd_{x_j} -\bd_{x'_j})(\bd_{y_k} -\bd_{y'_k})\Big\}L(x, y; x', y')\nonumber\\
=&\ \int dxdydx'dy'\ \delta(x-x')\delta(y-y') \sum_{jk} \bd_{jk}a(x-y) \nonumber\\
&\ \times \Big\{ (\bd_{y_j} -\bd_{x_j})(\bd_{y_k'} -\bd_{x'_k})+(\bd_{x_j} + \bd_{y'_j})( \bd_{x'_k}+ \bd_{y_k})\Big\}L(x, y; x', y'). \nonumber
\end{align*}
Writing in matrix notation  (with $\vect{A} =\vect{A}(x, y)$, and $\nabla$ a column vector)
\begin{align}
\frac{1}{2}\eqref{morawetz2-term3}=&\ \int dxdydx'dy'\ \delta(x-x')\delta(y-y')\nonumber \\
&\ \times \vect{A} : \Big\{(\grad_x-\grad_{y})(\grad_{x'}-\grad_{y'})^T L(x, y; x', y') \label{positive-term1}\\
&\ +(\grad_x\grad^T_{x'}+\grad_y \grad^T_{y'}) L(x, y; x', y')  \label{positive-term2}\\
&\ +(\grad_x\grad^T_{y}+\grad_{x'} \grad^T_{y'}) L(x, y; x', y')\Big\}  . \label{positive-term3}
\end{align}
\end{subequations}
Since $L$ is a positive operator, then it has a unique positive square root $\sqrt{L}$ such that $L = \sqrt{L}\circ \sqrt{L}$. In particular, we can now write
\begin{align*}
\eqref{positive-term1} =&\  \int dxdyd x_2 dy_2 dx'dy' dx_2' dy_2'\ \delta(x-x')\delta(y-y')\delta(x_2-x_2')\delta(y_2-y_2')\\
&\ \times \vect{A}:\left\{(\grad_x-\grad_{y}) \sqrt{L}(x, y; x'_2, y'_2)\overline{(\grad_{x'}-\grad_{y'})^T \sqrt{L}(x', y'; x_2, y_2)}\right\}\\
=&\  \int dxdyd x_2 dy_2 dx'dy' dx_2' dy_2'\ \delta(x-x')\delta(y-y')\delta(x_2-x_2')\delta(y_2-y_2')\\
&\ \times \overline{(\grad_{x'}-\grad_{y'})^T \sqrt{L}(x', y'; x_2, y_2)}\vect{A}(\grad_x-\grad_{y}) \sqrt{L}(x, y; x'_2, y'_2)\\
=&\ \llp{\vect{A}^\frac{1}{2}(\grad_y-\grad_x) \sqrt{L} }_\text{HS}^2\ge 0.
\end{align*}
The same argument holds for \eqref{positive-term2}, that is
\begin{align*}
\eqref{positive-term2} =  \llp{\vect{A}^\frac{1}{2}\grad_x \sqrt{L} }_\text{HS}^2+ \llp{\vect{A}^\frac{1}{2}\grad_y \sqrt{L} }_\text{HS}^2.
\end{align*}
For the final term, we need  the observation $\sqrt{L}(x, y; x', y') = \sqrt{\Gamma}(x; x')\sqrt{\Gamma}(y; y')$. Then it follows that
\begin{align*}
\eqref{positive-term3} =&\ \int dxdyd x_2 dy_2 dx'dy' dx_2' dy_2'\ \delta(x-x')\delta(y-y')\delta(x_2-x_2')\delta(y_2-y_2')\\
&\ \times \vect{A}:\Big\{\grad_x\sqrt{\Gamma}(x; x'_2)\grad^T_y\sqrt{\Gamma}(y; y'_2)\sqrt{\Gamma}(x_2; x')\sqrt{\Gamma}(y_2; y')\\
&\ +\sqrt{\Gamma}(x; x'_2)\sqrt{\Gamma}(y; y'_2)\grad_{x'}\sqrt{\Gamma}(x_2; x')\grad^T_{y'}\sqrt{\Gamma}(y_2; y')\Big\}\\
=&\ \int dxdyd x_2 dy_2 dx'dy' dx_2' dy_2'\ \delta(x-x')\delta(y-y')\delta(x_2-x_2')\delta(y_2-y_2')\\
&\ \times \Big\{\left(\grad_{x} \sqrt{L}(x, y_2; x_2', y')\right)^T\vect{A}\grad_{y} \sqrt{L}(x_2, y; x', y'_2)\\
&\ + \overline{\grad_{y'} \sqrt{L}(x_2', y'; x, y_2)}^T\vect{A}\overline{\grad_{x'} \sqrt{L}(x', y'_2; x_2, y)}\Big\}.
\end{align*}
Finally, by Cauchy-Schwarz inequality, we have that
\begin{align*}
|\eqref{positive-term3}| \ge - 2\llp{\vect{A}^\frac{1}{2}\grad_x \sqrt{L} }_\text{HS}\llp{\vect{A}^\frac{1}{2}\grad_y \sqrt{L} }_\text{HS}.
\end{align*}
Hence the desired result follows.
\end{proof}

\begin{proposition}
Let $\Gamma(t)$ be a  smooth global solution to \eqref{bbgky-gamma}
 with ${\rm tr} \Gamma(t)=1$ and $E(t) \le C$ (see \eqref{b5-nbr0}, \eqref{energy}).
 Then the following estimate
\begin{align}
\int dtdx\ |\Gamma(t,x, x)|^2 \lesssim 1
\end{align}
holds uniformly in $N$ and depends only on the initial data. Moreover, we also have the estimate
\begin{align}
\llp{\phi}_{L^4(dtdx)}\lesssim 1.
\end{align}
\end{proposition}

\begin{proof}
By the above lemmas, it immediately follows that
\begin{align}
8\pi \int^T_{-T}dt\int dx\ \rho(t, x)^2 \leq M^a(T)-M^a(-T).
\end{align}
To complete the argument, let us recall that $\Gamma(x; x')= \bar\phi(x)\phi(x') + N^{-1} (\overline{\sh(k)} \circ \sh(k))(x; x')$, then we see that
\begin{subequations}
\begin{align}
M^a(t) =& \int dxdy\  \rho(y)  \frac{x-y}{|x-y|} \cdot \Im\left(\bar \phi(x) \grad\phi(x)  \right)\\
&+\frac{1}{N}\int dxdy\  \rho(y)  \frac{x-y}{|x-y|} \cdot \Im\left(\overline{\sh(k)}\circ \grad\sh(x)  \right).
\end{align}
\end{subequations}
Finally, by a standard momentum-type estimate (see Lemma A.10 in \cite{Tao}, we see that
\begin{align*}
|M(t)| \le C \int dy\  \rho(t, y) \left\{\llp{|\grad|^{1/2}\phi(t)}_{L^2}^2+\frac{1}{N}\llp{|\grad_x|^{1/2}\sh(k_t)}^2_{L^2}\right\}.
\end{align*}
Finally, by the conservation of numbers and energy, we have the desired estimate.
\end{proof}

\section{Second proof of Proposition \ref{dispersive}\label{secondproof}}

Since this proposition is the main new technical ingredient of our paper, we give a second proof which is not based
on the kernel of the operator $A(x+2tD)$ (Weyl calculus), but rather on the Green's function.

We would like to show the following estimate,
\begin{align*}
\sup_{x_{1}}
\Big\Vert e^{it\big(\Delta_{x_{1}}+\Delta_{x_{2}}\big)}f(x_{1},x_{2})
\Big\Vert_{L^{2}(dx_{2})}\leq \frac{C}{t^{\frac{3}{2}}}
\big\Vert f\big\Vert_{L^{1}(dx_{1-2})L^{2}(dx_{1+2})}
\end{align*}
where (for convenience) we set
\begin{align*}
x_{1+2}:=\frac{x_{1}+x_{2}}{\sqrt{2}}\quad ,\quad
x_{1-2}:=\frac{x_{1}-x_{2}}{\sqrt{2}}.
\end{align*}
As in the first proof, we take the singular value decomposition of
$f(x_{1},x_{2})$ in the rotated $(x_{1-2},x_{1+2})$ variables
and write
\begin{align*}
f(x_{1},x_{2})=\sum_{k}u_{k}\left(\frac{x_{1}-x_{2}}{\sqrt{2}}\right)
v_{k}\left(
\frac{x_{1}+x_{2}}{\sqrt{2}}\right)
\end{align*}
where $\{v_{k}\}$ are orthonormal and $\{u_{k}\}$ are orthogonal.
The evolution equation can be written with the help of the
Green's functions as follows,
\begin{align*}
&e^{it\big(\Delta_{x_{1}}+\Delta_{x_{2}}\big)}f(x_{1},x_{2})
\\
&=\frac{1}{(4\pi t)^{3}}
\int\limits_{\RR^{3}\times\RR^{3}}dy_{1}dy_{2}
\sum_{k}\left\{u_{k}(y_{1})v_{k}(y_{2})
\exp\left(i\frac{\vert x_{1-2}-y_{1}\vert^{2}}{4t}
+i\frac{\vert x_{1+2}-y_{2}\vert^{2}}{4t}\right)
\right\}
\end{align*}
The phase in the exponential can be expanded,
\begin{align*}
\frac{\vert x_{1-2}-y_{1}\vert^{2}}{4t}+
\frac{\vert x_{1+2}-y_{2}\vert^{2}}{4t}
&=\frac{\vert x_{1}\vert^{2}+\vert x_{2}\vert^{2}+\vert y_{1}\vert^{2}
+\vert y_{2}\vert^{2}}{4t}
\\
&-\frac{x_{1}\cdot y_{1+2}}{2t}-\frac{x_{2}\cdot y_{2-1}}{2t}
\end{align*}
and in view of the above we redefine,
 \begin{align*}
& u_{k}(t,y_{1},x_{1}):=u_{k}(y_{1})
 \exp\left(i\frac{\vert y_{1}\vert^{2}-\sqrt{2}x_{1}\cdot y_{1}}{4t}\right)
 \\
 &v_{k}(t,y_{2},x_{1})
 :=v_{k}(y_{2})\exp\left(i
 \frac{\vert y_{2}\vert^{2}-\sqrt{2}x_{1}\cdot y_{2}}{4t}\right).
\end{align*}
Notice that
\begin{align*}
\big\{v_{k}(t,\cdot,x_{1})\big\}_{k}\quad {\rm is\ orthonormal.}
\end{align*}
Next we pick some function $A(x_{2})\in L^{2}(\RR^{3})$ and employ duality,
\begin{align*}
&\int_{\RR^{3}} dx_{2}
\left\{e^{it\big(\Delta_{x_{1}}+\Delta_{x_{2}}\big)}f(x_{1},x_{2})
A(x_{2})\right\}
\\
&=\frac{e^{i\frac{\vert x_{1}\vert^{2}}{4t}}}{(4\pi   t)^{3}}
\int\limits_{\RR^{3}\times\RR^{3}}
dy_{1}dy_{2}\sum_{k}
\left\{u_{k}(t,y_{1},x_{1})v_{k}(t,y_{2},x_{1})
\int\limits_{\RR^{3}}dx_{2}\left\{e^{-i\frac{x_{2}\cdot y_{2-1}}{2t}}
e^{i\frac{\vert x_{2}\vert^{2}}{4t}}A(x_{2})\right\}
\right\}
\\
&=\frac{e^{i\frac{\vert x_{1}\vert^{2}}{4t}}}{(4\pi   t)^{3}}
\int\limits_{\RR^{3}\times\RR^{3}}
dy_{1}dy_{2}\sum_{k}
\left\{u_{k}(t,y_{1},x_{1})v_{k}(t,y_{2},x_{1})
\widehat{A}\left(t,\frac{y_{2-1}}{2t}\right)
\right\}
\end{align*}
where we set,
\begin{align*}
&A(t,x_{2}):=e^{i\frac{\vert x_{2}\vert^{2}}{4t}}A(x_{2})
\\
&\widehat{A}(t,\xi)=\int_{\RR^{3}}dx_{2}\left\{
e^{-ix_{2}\cdot\xi}A(t,x_{2})\right\}.
\end{align*}
Let us now define
\begin{align*}
c_{k}(t,x_{1},y_{1}):=\int_{\RR^{3}}
dy_{2}\left\{v_{k}(t,y_{2},x_{1})\widehat{A}
\Big(t,\frac{y_{2}-y_{1}}{2\sqrt{2}t}
\right\}
\end{align*}
and the orthonormality of the set $\{v_{k}(t,\cdot,x_{1})\}$
imply
\begin{align*}
\sum_{k}\big\vert c_{k}(t,x_{1},y_{1})\big\vert^{2}
\leq C\int dy_{2}\left\{\Big\vert\widehat{A}\Big(
t,\frac{y_{2}-y_{1}}{2\sqrt{2}t}\Big)\Big\vert^{2}\right\}
=Ct^{3}\Vert A\Vert^{2}_{L^{2}(\RR^{3})}.
\end{align*}
Finally we have using Cauchy-Schwartz,
\begin{align*}
&\sup_{x_{1}\in\RR^{3}}
\left\vert \int_{\RR^{3}}dx_{2}\left\{
e^{i\big(\Delta_{x_{1}}+\Delta_{x_{2}}\big)}f(x_{1},x_{2})A(x_{2})
\right\}\right\vert
\\
&\leq \frac{C}{t^{3}}\int_{\RR^{3}}dy_{1}\left\{
\left(\sum_{k}\vert u_{k}(y_{1})\vert^{2}\right)^{\frac{1}{2}}
\left(\sum_{j}\vert c_{j}(t,x_{1},y_{1})\vert^{2}\right)^{\frac{1}{2}}
\right\}
\\
&\leq \frac{C}{t^{\frac{3}{2}}}\int_{\RR^{3}} dy_{1}\left(
\sum_{k}\vert u_{k}(y_{1})\vert^{2}\right)^{\frac{1}{2}}
\times \Vert A\Vert_{L^{2}(\RR^{3})}.
\end{align*}
The fact that $\{v_{k}\}$ are orthonormal imply that
\begin{align*}
\Vert f(x_{1},x_{2})\Vert_{L^{1}(dx_{1-2})L^{2}(dx_{1+2})}
=\left\Vert \sqrt{\sum_{k}\vert u_{k}(y_{1})\vert^{2}}
\right\Vert_{L^{1}(dy_{1})}.
\end{align*}

\section{Appendix: Proof of Theorem \ref{simple}}

 The detailed estimates for Term 1, Term 2 and Term 3 are slightly different (and irrelevant). They are
\begin{align*}
&\|\textmd{Term1}\|_{\mathcal{S}_{dual}^{\frac{8}{5}, \frac{4}{3}}[T_1, T_2]}
\lesssim  ||\Gamma(t,x,x)||_{L^4[T_1, T_2]L^2(dx)} \|\Lambda\|_{\mathcal{S}_{restrited}[T_1, T_2]}\\
&\|\nabla \textmd{Term1}\|_{\mathcal{S}_{dual}^{\frac{8}{5}, \frac{4}{3}}[T_1, T_2]}
\lesssim  N ||\Gamma(t,x,x)||_{L^4[T_1, T_2]L^2(dx)} \|\Lambda\|_{\mathcal{S}_{restrited}[T_1, T_2]}\\
&+||\Gamma(t,x,x)||_{L^4[T_1, T_2]L^2(dx)} \|\nabla \Lambda\|_{\mathcal{S}_{restrited}[T_1, T_2]}\\
&\|\nabla_x \nabla_y \textmd{Term1}\|_{\mathcal{S}_{dual}^{\frac{8}{5}, \frac{4}{3}}[T_1, T_2]}
\lesssim  N ||\Gamma(t,x,x)||_{L^4[T_1, T_2]L^2(dx)} \|\nabla \Lambda\|_{\mathcal{S}_{restrited}[T_1, T_2]}\\
&+||\Gamma(t,x,x)||^{\frac{1}{2}}_{L^2[T_1, T_2]L^2(dx)} \\
    &\cdot \left(||\nabla_x\nabla_y\Lambda(t,x,y)||_{L^{\frac{8}{3}}[T_1, T_2]L^4(dx)L^2(dy)}+
    ||\nabla_x\nabla_y\Lambda(t,x,y)||_{L^{\frac{8}{3}}[T_1, T_2]L^4(dy)L^2(dx)}\right).
\end{align*}
\begin{align*}
&\|\textmd{Term2}\|_{\mathcal{S}_{dual}^{\frac{8}{5}, \frac{4}{3}}[T_1, T_2]}
\lesssim N^{\frac{1}{2}} ||\Gamma(t,x,x)||^{\frac{1}{2}}_{L^4[T_1, T_2]L^2(dx)} \|\Lambda\|_{\mathcal{S}_{restrited}[T_1, T_2]}\\
&\|\nabla \textmd{Term2}\|_{\mathcal{S}_{dual}^{\frac{8}{5}, \frac{4}{3}}[T_1, T_2]}
\lesssim  N^{\frac{3}{2}} ||\Gamma(t,x,x)||^{\frac{1}{2}}_{L^4[T_1, T_2]L^2(dx)} \|\Lambda\|_{\mathcal{S}_{restrited}[T_1, T_2]}\\
&+N^{\frac{1}{2}}||\Gamma(t,x,x)||^{\frac{1}{2}}_{L^4[T_1, T_2]L^2(dx)} \|\nabla \Lambda\|_{\mathcal{S}_{restrited}[T_1, T_2]}\\
&\|\nabla_x \nabla_y \textmd{Term2}\|_{\mathcal{S}_{dual}^{\frac{8}{5}, \frac{4}{3}}[T_1, T_2]}
\lesssim  N^{\frac{3}{2}} ||\Gamma(t,x,x)||_{L^4[T_1, T_2]L^2(dx)} \|\nabla \Lambda\|_{\mathcal{S}_{restrited}[T_1, T_2]}.
\end{align*}
\begin{align*}
&\|\textmd{Term3}\|_{\mathcal{S}_{dual}^{\frac{8}{5}, \frac{4}{3}}[T_1, T_2]}
\lesssim  ||\Gamma(t,x,x)||_{L^4[T_1, T_2]L^2(dx)} \|\Lambda\|_{\mathcal{S}_{restrited}[T_1, T_2]}\\
&\|\nabla \textmd{Term3}\|_{\mathcal{S}_{dual}^{\frac{8}{5}, \frac{4}{3}}[T_1, T_2]}\\
 &\lesssim N||\Gamma(t, x,x)||_{L^4[T_1, T_2] L^2(dx)}||\Lambda||_{\mathcal{S}_{restrited}[T_1, T_2]}\\
 & +||\Gamma(t, x,x)||_{L^4[T_1, T_2] L^2(dx)} ||\nabla\Lambda||_{\mathcal{S}_{restrited}[T_1, T_2]}\\
  & + N^{\frac{3}{4}}||\Gamma(t, x,x)||^{\frac{1}{2}}_{L^2[T_1, T_2] L^2(dx)}||\Lambda||_{\mathcal{S}_{restrited}[T_1, T_2]}\\
&\|\nabla_x \nabla_y \textmd{Term3}\|_{\mathcal{S}_{dual}^{\frac{8}{5}, \frac{4}{3}}[T_1, T_2]}
\lesssim   N ||\Gamma(t,x,x)||_{L^4[T_1, T_2]L^2(dx)} \|\nabla \Lambda\|_{\mathcal{S}_{restrited}[T_1, T_2]}.
\end{align*}

To go from here to Theorem \ref{simple}, we estimate
\begin{align*}
 &||\Gamma(t,x,x)||_{L^4[T_1, T_2]L^2(dx)} \le ||\Gamma(t,x,x)||^{\frac{1}{2}}_{L^2[T_1, T_2]L^2(dx)}||\Gamma(t,x,x)||^{\frac{1}{2}}_{L^{\infty}[T_1, T_2]L^2(dx)}\\
 &\lesssim ||\Gamma(t,x,x)||^{\frac{1}{2}}_{L^2[T_1, T_2]L^2(dx)}.
 \end{align*}
 We present the detailed proofs, split into several propositions.

 The estimate for $\textmd{Term1}$
is an immediate consequence of H\"older's inequality, the Leibniz rule and $V_N(x) = N^{3 \beta} V(N^{\beta} x)$ with $\beta \le 1$.
 \begin{proposition} \label{term1} For any time interval $[T_1, T_2]$
 \begin{align*}
     &||\left(V_N \ast \Gamma(t,x,x)\right)\Lambda(t,x,y)||_{L^{\frac{8}{5}}[T_1, T_2]L^{\frac{4}{3}}(dx)L^2(dy)} \\
    &+||\left(V_N \ast \Gamma(t,y, y)\right)\Lambda(t,x,y)||_{L^{\frac{8}{5}}[T_1, T_2]L^{\frac{4}{3}}(dy)L^2(dx)} \\
    &\lesssim ||\Gamma(t,x,x)||_{L^4[T_1, T_2]L^2(dx)} \\
    &\cdot \left(||\Lambda(t,x,y)||_{L^{\frac{8}{3}}[T_1, T_2]L^4(dx)L^2(dy)}+
    ||\Lambda(t,x,y)||_{L^{\frac{8}{3}}[T_1, T_2]L^4(dy)L^2(dx)}\right)
    .
    \end{align*}
    while
     \begin{align*}
     &||\nabla_{x, y} \big(V_N \ast \Gamma(t,x,x)\Lambda(t,x,y)\big)||_{L^{\frac{8}{5}}[T_1, T_2]L^{\frac{4}{3}}(dx)L^2(dy)} \\
    &+||\nabla_{x, y}\big(V_N \ast \Gamma(t,y, y)\Lambda(t,x,y)\big)||_{L^{\frac{8}{5}}[T_1, T_2]L^{\frac{4}{3}}(dy)L^2(dx)} \\
    &\lesssim N ||\Gamma(t,x,x)||_{L^4[T_1, T_2]L^2(dx)} \\
    &\cdot \left(||\Lambda(t,x,y)||_{L^{\frac{8}{3}}[T_1, T_2]L^4(dx)L^2(dy)}+
    ||\Lambda(t,x,y)||_{L^{\frac{8}{3}}[T_1, T_2]L^4(dy)L^2(dx)}\right)\\
    &+ ||\Gamma(t,x,x)||_{L^4[T_1, T_2]L^2(dx)} \\
    &\cdot \left(||\nabla_{x, y}\Lambda(t,x,y)||_{L^{\frac{8}{3}}[T_1, T_2]L^4(dx)L^2(dy)}+
    ||\nabla_{x, y}\Lambda(t,x,y)||_{L^{\frac{8}{3}}[T_1, T_2]L^4(dy)L^2(dx)}\right)
    \end{align*}
    and
     \begin{align*}
     &||\nabla_x \nabla_y \big(V_N \ast \Gamma(t,x,x)\Lambda(t,x,y)\big)||_{L^{\frac{8}{5}}[T_1, T_2]L^{\frac{4}{3}}(dx)L^2(dy)} \\
    &+||\nabla_x \nabla_y\big(V_N \ast \Gamma(t,y, y)\Lambda(t,x,y)\big)||_{L^{\frac{8}{5}}[T_1, T_2]L^{\frac{4}{3}}(dy)L^2(dx)} \\
    &\lesssim N ||\Gamma(t,x,x)||_{L^4[T_1, T_2]L^2(dx)} \\
    &\cdot \left(||\nabla_{x, y}\Lambda(t,x,y)||_{L^{\frac{8}{3}}[T_1, T_2]L^4(dx)L^2(dy)}+
    ||\nabla_{x, y}\Lambda(t,x,y)||_{L^{\frac{8}{3}}[T_1, T_2]L^4(dy)L^2(dx)}\right)\\
    &+||\Gamma(t,x,x)||_{L^4[T_1, T_2]L^2(dx)} \\
    &\cdot \left(||\nabla_x\nabla_y\Lambda(t,x,y)||_{L^{\frac{8}{3}}[T_1, T_2]L^4(dx)L^2(dy)}+
    ||\nabla_x\nabla_y\Lambda(t,x,y)||_{L^{\frac{8}{3}}[T_1, T_2]L^4(dy)L^2(dx)}\right).
    \end{align*}

    \end{proposition}

 The propositions that follow are slightly more involved variants of the above argument.

In order to estimate $\textmd{Term2}$, we will use
\begin{proposition} \label{notroubleterm} For any time interval $[T_1, T_2]$,
\begin{align*}
&\big\|\left(V_N\Lambda \right)\circ \Gamma\big\|_{L^{\frac{8}{5}}([T_1, T_2])L^{\frac{4}{3}}(dx)L^2(dy)}\\
& \lesssim \|V_N\|^{\frac{1}{2}}_{L^{\frac{3}{2}}} \|\Lambda\|_{L^2[T_1, T_2]L^6(d(x-y)) L^2(d(x+y))}\|V_N\|^{\frac{1}{2}}_{L^1}
\|\Gamma(t, x, x)\|^{\frac{1}{2}}_{L^4[T_1, T_2]L^2}\|\Gamma(t, x, x)\|^{\frac{1}{2}}_{L^{\infty}L^1}\\
& \lesssim
 N^{\frac{1}{2}} \|\Gamma(t, x, x)\|^{\frac{1}{2}}_{L^4[T_1, T_2]L^2} \|\Lambda\|_{L^2[T_1, T_2]L^6(d(x-y)) L^2(d(x+y))}
\end{align*}
an also
\begin{align*}
&\big\|\bar \Gamma \circ \left(V_N\Lambda \right)\big\|_{L^{\frac{8}{5}}([T_1, T_2])L^{\frac{4}{3}}(dy)L^2(dx)}\\
& \lesssim
 N^{\frac{1}{2}} \|\Gamma(t, x, x)\|^{\frac{1}{2}}_{L^4[T_1, T_2]L^2} \|\Lambda\|_{L^2[T_1, T_2]L^6(d(x-y)) L^2(d(x+y))}.
\end{align*}
\end{proposition}
\begin{proof}
We have the pointwise estimate
\begin{align}
&\big|\left(V_N\Lambda \circ \Gamma\right)(t, x, y)\big|= \big|\int V_N(x-z) \Lambda (t, x, z)\Gamma(t, z, y) dz \big|\notag\\
& \le \left(\int V_N(x-z) |\Lambda (t, x, z)|^2 dz \right)^{\frac{1}{2}}
\left(\int V_N(x-z) |\Gamma (t, z, y)|^2 dz \right)^{\frac{1}{2}} \label{badterm}\\
&:=A(t, x) B(t, x, y).\notag
\end{align}
Thus
\begin{align*}
\big\|\left(V_N\Lambda \right)\circ \Gamma\big\|_{L^{\frac{8}{5}}([T_1, T_2])L^{\frac{4}{3}}(dx)L^2(dy)}\le
\|A\|_{L^2[T_1, T_2]L^2}\|B\|_{L^8L^4L^2}
\end{align*}
and
\begin{align*}
\|A\|_{L^2[T_1, T_2]L^2} \le \|V_N\|^{\frac{1}{2}}_{L^{\frac{3}{2}}} \|\Lambda\|_{L^2[T_1, T_2]L^6(d(x-y)) L^2(d(x+y))}.
\end{align*}
Also, using \eqref{gammaprop},
 \begin{align*}
& \bigg\|\left(\int V_N(x-z) |\Gamma (t, z, y)|^2 dz \right)^{\frac{1}{2}}\bigg\|_{L^{8 }[T_1, T_2] L^4(dx) L^2(dy)}\\
&\le
\bigg\|\left(\int V_N(x-z) |\Gamma (t, z, z)| dz |\Gamma(t, y, y)|\right)^{\frac{1}{2}}\bigg\|_{L^{8 }[T_1, T_2] L^4(dx) L^2(dy)}\\
&\le
 \|V_N\|^{\frac{1}{2}}_{L^1}
\|\Gamma(t, x, x)\|^{\frac{1}{2}}_{L^4[T_1, T_2]L^2}\|\Gamma(t, x, x)\|^{\frac{1}{2}}_{L^{\infty}L^1}.
 \end{align*}
The proof of the second estimate is similar.

\end{proof}

Next, we need the above estimate with derivatives.

\begin{proposition} \label{derivtroubleterm} For any time interval $[T_1, T_2]$,
\begin{align*}
&\big\|\left(V_N\Lambda \right)\circ\nabla_y  \Gamma\big\|_{L^{\frac{8}{5}}([T_1, T_2])L^{\frac{4}{3}}(dx)L^2(dy)}\\
& \lesssim \|V_N\|^{\frac{1}{2}}_{L^{\frac{3}{2}}} \|\Lambda\|_{L^2[T_1, T_2]L^6(d(x-y)) L^2(d(x+y))}\|V_N\|^{\frac{1}{2}}_{L^1}
\|\Gamma(t, x, x)\|^{\frac{1}{2}}_{L^4[T_1, T_2]L^2}\|E_k\|^{\frac{1}{2}}_{L^{\infty}L^1}\\
& \lesssim
 N^{\frac{1}{2}} \|\Gamma(t, x, x)\|^{\frac{1}{2}}_{L^4[T_1, T_2]L^2} \|\Lambda\|_{L^2[T_1, T_2]L^6(d(x-y)) L^2(d(x+y))}
\end{align*}
(we used Proposition \ref{energy1}).
Thus, using the Leibniz rule,
\begin{align*}
&\big\|\nabla_{x, y}\left(\left(V_N\Lambda \right)\circ \Gamma\right)\big\|_{L^{\frac{8}{5}}([T_1, T_2])L^{\frac{4}{3}}(dx)L^2(dy)}\\
& \lesssim
 N^{\frac{3}{2}} \|\Gamma(t, x, x)\|^{\frac{1}{2}}_{L^4[T_1, T_2]L^2} \|\Lambda\|_{L^2[T_1, T_2]L^6(d(x-y)) L^2(d(x+y))}\\
 &+  N^{\frac{1}{2}} \|\Gamma(t, x, x)\|^{\frac{1}{2}}_{L^4[T_1, T_2]L^2} \|\nabla_{x, y}\Lambda\|_{L^2[T_1, T_2]L^6(d(x-y)) L^2(d(x+y))}
 \end{align*}
 and
 \begin{align*}
&\big\|\nabla_x \nabla_y\left(\left(V_N\Lambda \right)\circ \Gamma\right)\big\|_{L^{\frac{8}{5}}([T_1, T_2])L^{\frac{4}{3}}(dx)L^2(dy)}\\
& \lesssim
 N^{\frac{3}{2}} \|\Gamma(t, x, x)\|^{\frac{1}{2}}_{L^4[T_1, T_2]L^2} \|\nabla_{x, y}\Lambda\|_{L^2[T_1, T_2]L^6(d(x-y)) L^2(d(x+y))}.
 \end{align*}
A similar estimate holds for
 \begin{align*}
&\big\|\nabla_{x, y} \left(\bar \Gamma \circ \left(V_N\Lambda \right)\right)\big\|_{L^{\frac{8}{5}}([T_1, T_2])L^{\frac{4}{3}}(dy)L^2(dx)}.
\end{align*}
\end{proposition}
\begin{proof}
The argument is similar to the previous proof, with minor modifications. We have the pointwise estimate
\begin{align}
&\big|\left(V_N\Lambda \circ \nabla_y \Gamma\right)(t, x, y)\big|= \big|\int V_N(x-z) \Lambda (t, x, z) \nabla_y \Gamma(t, z, y) dz \big|\notag\\
& \le \left(\int V_N(x-z) |\Lambda (t, x, z)|^2 dz \right)^{\frac{1}{2}}
\left(\int V_N(x-z) |\nabla_y\Gamma (t, z, y)|^2 dz \right)^{\frac{1}{2}} \label{badterm}\\
&:=A(t, x) C(t, x, y)\notag
\end{align}
and
\begin{align*}
\big\|\left(V_N\Lambda \right)\circ \nabla_y \Gamma\big\|_{L^{\frac{8}{5}}([T_1, T_2])L^{\frac{4}{3}}(dx)L^2(dy)}\le
\|A\|_{L^2[T_1, T_2]L^2}\|C\|_{L^8L^4L^2}.
\end{align*}
For $A$, we have already noticed
\begin{align*}
\|A\|_{L^2[T_1, T_2]L^2} \le \|V_N\|^{\frac{1}{2}}_{L^{\frac{3}{2}}} \|\Lambda\|_{L^2[T_1, T_2]L^6(d(x-y)) L^2(d(x+y))}.
\end{align*}
For $C$, we use \eqref{gammaprop2}:
 \begin{align*}
& \bigg\|\left(\int V_N(x-z) |\nabla_y \Gamma (t, z, y)|^2 dz \right)^{\frac{1}{2}}\bigg\|_{L^{8 }[T_1, T_2] L^4(dx) L^2(dy)}\\
&\le
\bigg\|\left(\int V_N(x-z) |\Gamma (t, z, z)| dz E_k(t, y)\right)^{\frac{1}{2}}\bigg\|_{L^{8 }[T_1, T_2] L^4(dx) L^2(dy)}\\
&\le
 \|V_N\|^{\frac{1}{2}}_{L^1}
\|\Gamma(t, x, x)\|^{\frac{1}{2}}_{L^4[T_1, T_2]L^2}\|E_k\|^{\frac{1}{2}}_{L^{\infty}L^1}.
 \end{align*}

\end{proof}

Next, we discuss $\textmd{Term3}$.
 \begin{proposition} \label{vGammaLambda}  For any time interval $[T_1, T_2]$
 \begin{align*}
  &||\int (V_N (x-z)\bar \Gamma)(x,z) \Lambda(z,y) dz||_{L^{\frac{8}{5}}[T_1, T_2]L^{\frac{4}{3}}(dx)L^2(dy)}\\
    &\lesssim ||\Gamma(t, x,x)||_{L^4[T_1, T_2] L^2(dx)} ||\Lambda||_{L^{\frac{8}{3}}[T_1, T_2]L^4(dx)L^2(dy)}.
\end{align*}
and also
\begin{align*}
&  ||\int  \Lambda(x,z)(V_N(z-y)\Gamma)(z,y) dz||_{L^{\frac{8}{5}}[T_1, T_2]L^{\frac{4}{3}}(dx)L^2(dy)}\\
& \lesssim ||\Gamma(t, y,y)||_{L^4[T_1, T_2] L^2(dy)} ||\Lambda||_{L^{\frac{8}{3}}[T_1, T_2]L^4(dy)L^2(dx)}.
  \end{align*}
\end{proposition}
 \begin{proof}
Using \ref{gammaprop}
together with H{\"o}lder's inequality and Young's inequality, we have
\begin{equation}
  ||\int (V_N \bar\Gamma)(x,z) \psi(z) dz||_{L_x^{\frac{4}{3}}}  \lesssim ||\Gamma(x,x)||_{L^2} ||\psi||_{L^4}.
\end{equation}
Thus,  at fixed time, using $\psi(x)= \|\Lambda(x, \cdot)\|_{L^2(dy)}$,
\begin{equation}
  ||\int (V_N\Gamma)(x,z) \Lambda(z,y) dz||_{L^{\frac{4}{3}}(dx)L^2(dy)}  \lesssim ||\Gamma(x,x)||_{L^2} ||\Lambda||_{L^4L^2}.
\end{equation}
The proof is finished by using H\"older's inequality.
The argument for the second estimate is similar.

\end{proof}

Next, we introduce derivatives:

 \begin{proposition} \label{vGammaLambda'}
 \begin{align*}
 & ||\int (V_N(x-z) \nabla_x \bar \Gamma)(x,z) \Lambda(z,y) dz||_{L^{\frac{8}{5}}[T_1, T_2]L^{\frac{4}{3}}(dx)L^2(dy)} \\
  &\leq \|V_N\|_{L^{\frac{4}{3}}} ||\Gamma(t, x,x)||^{\frac{1}{2}}_{L^2[T_1, T_2] L^2(dx)}
  \|E_k\|^{\frac{1}{2}}_{L^\infty(dt) L^1(dx)} ||\Lambda||_{L^{\frac{8}{3}}[T_1, T_2]L^4(dx)L^2(dy)}\\
  & \lesssim N^{\frac{3}{4}}||\Gamma(t, x,x)||^{\frac{1}{2}}_{L^2[T_1, T_2] L^2(dx)}||\Lambda||_{L^{\frac{8}{3}}[T_1, T_2]L^4(dx)L^2(dy)}.
\end{align*}
where $E_k$ is the kinetic energy density, see \eqref{kin} and the estimate of Proposition \ref{energy1}.
Thus
\begin{align*}
 & ||\nabla_{x, y}\int V_N(x-z)  \bar \Gamma (x,z) \Lambda(z,y) dz||_{L^{\frac{8}{5}}[T_1, T_2]L^{\frac{4}{3}}(dx)L^2(dy)} \\
 &\lesssim N||\Gamma(t, x,x)||_{L^4[T_1, T_2] L^2(dx)}||\Lambda||_{L^{\frac{8}{3}}[T_1, T_2]L^4(dx)L^2(dy)}\\
 & +||\Gamma(t, x,x)||_{L^4[T_1, T_2] L^2(dx)} ||\nabla_y \Lambda||_{L^{\frac{8}{3}}[T_1, T_2]L^4(dx)L^2(dy)}\\
  & + N^{\frac{3}{4}}||\Gamma(t, x,x)||^{\frac{1}{2}}_{L^2[T_1, T_2] L^2(dx)}||\Lambda||_{L^{\frac{8}{3}}[T_1, T_2]L^4(dx)L^2(dy)}
 .
\end{align*}
and
\begin{align*}
 & ||\nabla_x \nabla_y\int V_N(x-z)  \bar \Gamma (x,z) \Lambda(z,y) dz||_{L^{\frac{8}{5}}[T_1, T_2]L^{\frac{4}{3}}(dx)L^2(dy)} \\
 &\lesssim N||\Gamma(t, x,x)||_{L^4[T_1, T_2] L^2(dx)}||\nabla_y\Lambda||_{L^{\frac{8}{3}}[T_1, T_2]L^4(dx)L^2(dy)}\\
  & + N^{\frac{3}{4}}||\Gamma(t, x,x)||^{\frac{1}{2}}_{L^2[T_1, T_2] L^2(dx)}||\nabla_y\Lambda||_{L^{\frac{8}{3}}[T_1, T_2]L^4(dx)L^2(dy)}
  .
\end{align*}

 Similar estimates hold for
$
  \int  \Lambda(x,z)\ \left(V_N(z-y)\Gamma(z,y)\right) dz$.
   \end{proposition}
 \begin{proof}
Using \eqref{gammaprop2} and arguing as in the previous proof, with
 $\psi(x)=\|\Lambda(x, \cdot)\|_{L^2(dy)}$,
we have,
\begin{align*}
&  ||\int (V_N \nabla_x \bar\Gamma)(x,z) \psi(z) dz||_{L^{\frac{4}{3}}}  \lesssim ||E_k^{\frac{1}{2}} V_N * \left(\Gamma(z, z)^{\frac{1}{2}} \psi(z)\right)||_{L^{\frac{4}{3}}}\\
& \le \|E_k^{\frac{1}{2}}\|_{L^2}
\|V_N * \left(\Gamma(z, z)^{\frac{1}{2}} \psi(z)\right)\|_{L^4}\\
& \le  \|E_k^{\frac{1}{2}}\|_{L^2}
\|V_N\|_{L^{\frac{4}{3}}}\|(\Gamma(z, z)^{\frac{1}{2}} \psi(z))\|_{L^2}\\
&\le \|V_N\|_{L^{\frac{4}{3}}} ||\Gamma( x,x)^{\frac{1}{2}}||_{ L^4(dx)}
  \|E_k\|^{\frac{1}{2}}_{ L^1(dx)} ||\Lambda||_{L^4(dx)L^2(dy)}.
\end{align*}
Now the result follows using H\"older's inequality in time.
The proof of the second estimate is similar.

\end{proof}

Finally, we need estimates for $(V_N\ast|\phi|^2)(x)\phi(x)\phi(y)$.

\begin{proposition}
\begin{align*}
 & \big\|(V_N\ast|\phi|^2)(x)\phi(x)\phi(y)\big\|_{L^{2}(dt)L^{\frac{6}{5}}(dx)L^2(dy)}
     +\big\|(V_N\ast|\phi|^2)(y)\phi(x)\phi(y)\big\|_{L^{2}(dt)L^{\frac{6}{5}}(dy)L^2(dx)}
     \\
    &\lesssim ||(V_N\ast|\phi|^2)(x)||_{L^2(dt)L^{2}(dx)}
    ||\phi||_{L^{\infty}(dt)L^{3}(dx)}||\phi||_{L^{\infty}(dt)L^2(dy)}\\
    &\lesssim 1\\
    & \big\|\nabla_{x, y}\left((V_N\ast|\phi|^2)(x)\phi(x) \phi(y)\right)\big\|_{L^{2}(dt)L^{\frac{6}{5}}(dx)L^2(dy)}\\
    &
     +\big\|\nabla_{x, y}\left((V_N\ast|\phi|^2)(y)\phi(x) \phi(y)\right)\big\|_{L^{2}(dt)L^{\frac{6}{5}}(dy)L^2(dx)}
 \\
    &\lesssim N\\
    & \big\|\nabla_x \nabla_y\left((V_N\ast|\phi|^2)(x)\phi(x) \phi(y)\right)\big\|_{L^{2}(dt)L^{\frac{6}{5}}(dx)L^2(dy)}\\
    &
     +\big\|\nabla_x \nabla_y\left((V_N\ast|\phi|^2)(y)\phi(x) \phi(y)\right)\big\|_{L^{2}(dt)L^{\frac{6}{5}}(dy)L^2(dx)}
 \\
    &\lesssim N.
    \end{align*}

\end{proposition}

\begin{proof}
    All the above can be proved using \eqref{b5-nbr0}, \eqref{energy} and \eqref{morawetz}.

Since $ \big\|\nabla_x \nabla_y  \Lambda \big\|_{L^{\infty}(dt)L^2(dx dy)} \lesssim \big\|\nabla_x \nabla_y  \Lambda \big\|_{{\mathcal S}_{restricted}}$, the proof of Theorem \ref{D^2} is complete.

\end{proof}


\begin{thebibliography}{10}

\bibitem{BBCFS}
V.~Bach, S.~Breteaux, T.~Chen, J.~Fr{\"o}hlich, and I.~M. Sigal, \emph{The
  time-dependent {Hartree-Fock-Bogoliubov} equations for bosons}, arXiv
  preprint arXiv:1602.05171 (2016), p. 1--46.

\bibitem{BSS}
N.~Benedikter, J.~Sok, and J.~P. Solovej, \emph{The {D}irac--{F}renkel
  principle for reduced density matrices, and the {B}ogoliubov--de {G}ennes
  equations}, Annales Henri Poincar{\'e} \textbf{19} (2018), no.~4, 1167--1214.

\bibitem{BCS}
C.~Boccato, S.~Cenatiempo, and B.~Schlein, \emph{Quantum many-body fluctuations
  around nonlinear {S}chr{\"o}dinger dynamics}, Annales Henri Poincar{\'e}
  \textbf{18} (2017), no.~1, 113--191.

\bibitem{aBPPS}
L.~Bo{\ss}mann, N.~Pavlovi{\'c}, P.~Pickl, and A.~Soffer, \emph{Higher order
  corrections to the mean-field description of the dynamics of interacting
  bosons}, Journal of Statistical Physics \textbf{178} (2020), no.~6,
  1362--1396.

\bibitem{BPPS}
L.~Bo{\ss}mann, S.~Petrat, P.~Pickl, and A.~Soffer, \emph{Beyond bogoliubov
  dynamics}, arXiv preprint arXiv:1912.11004 (2019), pp 1--62.

\bibitem{Bouclet}
J.~Bouclet and H.~Mizutani, \emph{Uniform resolvent and {S}trichartz estimates
  for {S}chr{\"o}dinger equations with critical singularities}, Transactions of
  the American Mathematical Society \textbf{370} (2018), no.~10, 7293--7333.

\bibitem{Bourgain98}
J.~Bourgain, \emph{Scattering in the energy space and below for 3{D} {NLS}},
  Journal d'Analyse Math{\'e}matique \textbf{75} (1998), no.~1, 267--297.

\bibitem{BNNS}
C.~Brennecke, P.~T. Nam, M.~Napi{\'o}rkowski, and B.~Schlein,
  \emph{Fluctuations of {N}-particle quantum dynamics around the nonlinear
  {S}chr{\"o}dinger equation}, Annales de l'Institut Henri Poincar{\'e} C,
  Analyse nonlin{\'e}aire \textbf{36} (2019), no.~5, 1201--1235.

\bibitem{CPTz}
T.~Chen, N.~Pavlovic, and N.~Tzirakis, \emph{Multilinear {M}orawetz identities
  for the {G}ross-{P}itaevskii hierarchy}, Contemp. Math \textbf{581} (2012),
  39--62.

\bibitem{Jackythesis}
J.~J.~W. Chong, \emph{Application of dispersive {PDE} techniques to the studies
  of the time-dependent {H}artree-{F}ock-{B}ogoliubov system for bosons}, Ph.D.
  thesis, University of Maryland, College Park, 2019.

\bibitem{jacky2}
\bysame, \emph{Dynamical {H}artree-{F}ock-{B}ogoliubov approximation of
  interacting bosons}, arXiv preprint arXiv:1711.00610v2 (2019), pp 1--56.

\bibitem{CKSTT}
J.~Colliander, M.~Keel, G.~Staffilani, H.~Takaoka, and T.~Tao, \emph{Existence
  globale et diffusion pour l'{\'e}quation de {S}chr{\"o}dinger nonlin{\'e}aire
  r{\'e}pulsive cubique sur $\mathbb{R}^3$ en dessous l'espace d'{\'e}nergie},
  Journ{\'e}es {\'e}quations aux d{\'e}riv{\'e}es partielles (2002), no.~10,
  1--15.

\bibitem{Dancona}
P.~D'Ancona, \emph{Kato smoothing and {S}trichartz estimates for wave equations
  with magnetic potentials}, Communications in Mathematical Physics
  \textbf{335} (2015), no.~1, 1--16.

\bibitem{FLLS}
R.~L. Frank, M.~Lewin, E.~H. Lieb, and R.~Seiringer, \emph{Strichartz
  inequality for orthonormal functions}, Journal of the European Mathematical
  Society \textbf{16} (2014), no.~7, 1507--1526.

\bibitem{GM2}
M.~Grillakis and M.~Machedon, \emph{Beyond mean field: On the role of pair
  excitations in the evolution of condensates}, Journal of Fixed Point Theory
  and Applications \textbf{14} (2013), no.~1, 91--111.

\bibitem{G-M2017}
\bysame, \emph{Pair excitations and the mean field approximation of interacting
  {Bosons}, {II}}, Communications in Partial Differential Equations \textbf{42}
  (2017), no.~1, 24--67.

\bibitem{G-M2019}
\bysame, \emph{Uniform in {N} estimates for a {B}osonic system of
  {H}artree--{F}ock--{B}ogoliubov type}, Communications in Partial Differential
  Equations \textbf{44} (2019), no.~12, 1431--1465.

\bibitem{GMM1}
M.~Grillakis, M.~Machedon, and D.~Margetis, \emph{Second-order corrections to
  mean field evolution of weakly interacting {Bosons}. {I}.}, Communications in
  Mathematical Physics \textbf{294} (2010), no.~1, 273--301.

\bibitem{JSS}
J-L Journ{\'e}, A.~Soffer, and C.~D. Sogge, \emph{Decay estimates for
  {S}chr{\"o}dinger operators}, Communications on Pure and Applied mathematics
  \textbf{44} (1991), no.~5, 573--604.

\bibitem{Kato}
T.~Kato, \emph{Non-existence of bound states with positive energy}, Journal of
  the Physical Society of Japan \textbf{14} (1959), no.~3, 382--382.

\bibitem{K-T}
M.~Keel and T.~Tao, \emph{Endpoint strichartz estimates}, American Journal of
  Mathematics \textbf{120} (1998), no.~5, 955--980.

\bibitem{Mizutani}
H.~Mizutani, J.~Zhang, and J.~Zheng, \emph{Uniform resolvent estimates for
  {S}chr{\"o}dinger operator with an inverse-square potential}, Journal of
  Functional Analysis \textbf{278} (2020), no.~4, 108350.

\bibitem{Rellich}
F.~Rellich, \emph{{\"U}ber das asymptotische {V}erhalten der {L}{\"o}sungen von
  {$\Delta u+ \lambda u= 0$} in unendlichen gebieten.}, Jahresbericht der
  Deutschen Mathematiker-Vereinigung \textbf{53} (1943), 57--65.

\bibitem{Rod-S}
I.~Rodnianski and B.~Schlein, \emph{Quantum fluctuations and rate of
  convergence towards mean field dynamics}, Communications in Mathematical
  Physics \textbf{291} (2009), no.~1, 31--61.

\bibitem{Tao}
T.~Tao, \emph{Local and global analysis of nonlinear dispersive and wave
  equations}, CBMS Regional Conference Series in Mathematics, no. 106, American
  Mathematical Society, 2006.

\bibitem{yajima}
K.~Yajima, \emph{The ${W}^{k, p}$ continuity of wave operators for
  {S}chr{\"o}dinger operators}, Journal of the Mathematical Society of Japan
  \textbf{47} (1995), no.~3, 551--581.

\end{thebibliography}
\end{document}